\documentclass[12pt]{article}
\usepackage{amsfonts,amssymb,amsmath,theorem}
\renewcommand{\thesubsection}{\thesection\kern0.08em\alph{subsection}}
\newtheorem{Pa}{Paper}[section]
\newtheorem{Tm}[Pa]{{\bf Theorem}}
\newtheorem{La}[Pa]{{\bf Lemma}}
\newtheorem{Cy}[Pa]{{\bf Corollary}}
\newtheorem{Pn}[Pa]{{\bf Proposition}}

{
\theorembodyfont{\normalfont}
\newtheorem{I}[Pa]{{\bf}}
\newtheorem{Rk}[Pa]{{\bf Remark}}
\newtheorem{Ee}[Pa]{{\bf Example}}
\newtheorem{Dn}[Pa]{{\bf Definition}}

}

\newcommand\norm[1]{\left\|#1\right\|}           
\newcommand\qed{\ifhmode\unskip\nobreak\fi\quad  
   \ifmmode\square\else\hbox{$\square$}\fi}      
\let\emptyset=\varnothing
\let\cal=\mathcal

\newcommand\xrlarrow[3][]{\genfrac{}{}{0pt}{3}%
{\displaystyle\smash[b]{\xrightarrow[#1\phantom{#3}#1]{#2}}}%
{\displaystyle\smash[t]{\xleftarrow[#1{#3}#1]{\phantom{#2}}}}}%

\numberwithin{equation}{section}

\begin{document}
\begin{center}

{\bfseries \large Crossed product of a  $C^*$-algebra  by an
endomorphism,\\[3pt]
coefficient algebras and transfer operators}

\bigbreak
A.B. Antonevich,  V.I. Bakhtin, A.V. Lebedev\\

\bigskip
Belarus State University / University of Bialystok\\
Belarus State University;\\
Belarus State University / University of Bialystok
\end{center}

\vspace{5mm}
\quad \parbox{14,5cm}{\small \hspace{0.5cm}
The paper  presents  a  construction  of  the  crossed  product  of  a
$C^*$-algebra\\ by an endomorphism generated by partial isometry.


 }

\vspace{5mm}
\tableofcontents

\section{Introduction}
\label{Intro}

Given a $C^*$-algebra $A$ and an endomorphism $\alpha$ there is a
number of ways to construct a new $C^*$-algebra (an extension of $A$)
called the crossed product. Among the successful constructions of this
sort one should mention, for example, the constructions developed by
J.\ Cuntz and W.\ Krieger \cite{cuntz,CK}, W.\,L.\ Paschke
\cite{Paschke}, P.\,J.\ Stacey \cite{Stacey}, G.\,J.\ Murphy
\cite{Murphy}, R.\ Exel \cite{exel}, and B.\,K.\ Kwasniewski
\cite{kwa}. Here the Exel's crossed product \cite{exel} is the most
general one since all the others can be reduced to it by means of this
or that procedure. At the same time the foregoing statement is 'not
completely true'. By saying this we mean the following. The Exel's
construction (see \ref{Ex}) starts with the {\em two} objects --- an
{\em endomorphism} $\alpha$ of a unital $C^*$-algebra $A$ and a {\em
transfer operator} $\cal L$ satisfying the prescribed relations. The
crossed product is then defined as the enveloping universal
$C^*$-algebra generated by $A$ and the operator $S$ that generates the
transfer operator and satisfies the given relations. If one considers
the Kwasniewski's crossed product (see Subsection \ref{Kwa}) then
there is {\em no} transfer operator among the starting objects. The
same is true for the Cuntz-Krieger algebra (see Subsection \ref{C-K}).
The recent paper by R.\ Exel \cite{exel-inter} makes the situation even
more intriguing (in this paper the {\em endomorphism disappears}). The
Exel's crossed product also possesses a certain 'drawback' --- as is
shown by N.\ Brownlowe and I.\ Raeburn \cite{Brow-Rae} it does not
always contain the initial algebra $A$ and therefore it is {\em not}
always an {\em extension} of $A$. All this means that the natural
question: 'What is the 'construction' of the crossed product in
general?' is still waiting for the answer.

In this article we investigate this question from one more side.

In the paper by A.\,V.\ Lebedev and A.\ Odzijewicz \cite{Leb-Odz} the
notion of the so-called {\em coefficient algebra} was introduced, and
it was shown that this object plays a principal role in the extensions
of $C^*$-algebras by partial isometries. In the paper by
V.\,I.\ Bakhtin and A.\,V.\ Lebedev \cite{Bakht-Leb} the criterion for a
$C^*$-algebra to be a coefficient algebra associated with a given
endomorphism was obtained. On the base of the results of these papers
one naturally arrives at the construction of a certain crossed
product. This construction is simpler and 'more natural' than the
Exel's one. The discussion of this crossed product is the main theme
of the present article.

The paper is organized as follows.

In Section \ref{crossed product} we recall the necessary results on
the coefficient algebras and transfer operators borrowed from
\cite{Leb-Odz}, \cite{exel} and \cite{Bakht-Leb} and introduce the
notion of the crossed product. Section \ref{faithful} is devoted to
the description of the internal structure of the crossed product
introduced. Here, in particular, we give a criterion for a
representation to be a faithful representation of the crossed product
and present the regular representation. In Section \ref{v-crossed
product} we compare the crossed product introduced in this article
with the already existing successful constructions. It appears that it
covers 'almost all' of them, in particular, in the 'most popular
situation' when all the powers of the transfer operator are generated
by partial isometries it coincides with the Exel's crossed product but
with {\em different} algebra, {\em different} endomorphism and {\em
different} transfer operator. The final Section \ref{crossed product
etc} is devoted to the discussion of the 'general crossed product
construction'. Though this section does not contain the 'exact'
results we consider it as being among the important parts of this
article since it contains the 'general' crossed product 'philosophy'
(from our point of view).

\section{Coefficient algebras,  transfer operators,\hfil\null\ \hbox{crossed product}\hfill}
\label{crossed product}

For the sake of completeness of the presentation we start with
recalling some definitions and facts concerning transfer operators and
coefficient algebras. The corresponding material is borrowed from
\cite{exel}, \cite{Leb-Odz} and \cite{Bakht-Leb}.

Let  $\cal A$ be a $C^*$-algebra with an identity $1$ and
$\delta\!:\cal A\to \cal A$ be an endomorphism of this algebra.
A linear map $\delta_*\!:\cal A\to \cal A$ is called a
{\em transfer operator\/} for the pair  $(\cal A,\delta)$ if it is
continuous and  positive and such that
\begin{equation}\label{b,,2}
\delta_*(\delta(a)b) =a\delta_*(b),\qquad a,b\in\cal A.
\end{equation}

\begin{Pn}
\label{Ex2.3} {\upshape \cite[Proposition 2.3]{exel} }\ \
Let $\delta_*$ be a transfer operator for the pair $(\cal A,\delta)$.
Then the following are equivalent:

\smallskip
\quad\ \llap{$(i)$}\ \
the composition $E = \delta \circ \delta_*$ is a conditional
expectation onto $\delta (\cal A)$,

\smallskip
\quad\ \llap{$(ii)$}\ \ $\delta \circ \delta_* \circ \delta = \delta$,

\smallskip
\quad\ \llap{$(iii)$}\ \ $\delta (\delta_*(1)) = \delta (1)$.
\end{Pn}

If the equivalent conditions of Proposition  \ref{Ex2.3} hold then
R.\ Exel calls $\delta_*$ a {\em non-degenerate} transfer operator.

The transfer operator $\delta_*$  is called {\em complete,} if
\begin{equation}
\label{b,,3}
\delta\delta_*(a) =\delta(1)a\delta(1),\qquad a\in\cal A.
\end{equation}
Observe that a complete transfer operator is non-degenerate. Indeed,
\eqref{b,,3} implies
\[
\delta\delta_*\delta(a) = \delta(1)\delta (a) \delta(1) = \delta (a)
\]
and so condition (ii) of Proposition  \ref{Ex2.3} is satisfied.

The next result presents the criteria for the existence of a complete
transfer operator.

\begin{Tm}\label{complete}
{\upshape \cite[Theorem  2.8]{Bakht-Leb}}\ \
Let $\cal A$ be a $C^*$-algebra with an identity $1$ and $\delta \!:
{\cal A} \to {\cal A}$ be an endomorphism of $\cal A$. The following
are equivalent:

\medskip
\quad\llap{$1)$}\ \
there exists a complete transfer operator $\delta_*$ $($for
$({\cal A}, \delta )$$),$

\medskip
\quad\llap{$2)$}\qquad\llap{$(i)$}\ \,
there exists a non-degenerate transfer operator $\delta_*$ and

\smallskip
\quad\qquad\llap{$(ii)$}\ \,
$\delta ({\cal A})$ is a hereditary subalgebra of $\cal A$$;$

\medskip
\quad\llap{$3)$}\qquad\llap{$(i)$}\ \,
there exists a central orthogonal projection  $P\in \cal A$ such that

\smallskip
\qquad\qquad\quad\llap{$a)$}\ $\delta (P) = \delta (1)$,

\smallskip
\qquad\qquad\quad\llap{$b)$}\ the mapping $\delta \!: P {\cal A} \to
\delta ({\cal A})$ is a $^*$-isomorphism, and

\smallskip
\quad\qquad\llap{$(ii)$}\ \, $\delta ({\cal A}) = \delta (1){\cal A}
\delta (1)$.

\medbreak\noindent
Moreover the objects in $1)$ -- $3)$ are defined in a unique way
$($i.\,e.\ the transfer operator $\delta_*$ in $1)$ and $2)$ is
unique and the projection $P$ in $3)$ is unique as well$)$ and
\begin{equation}
\label{P} P=\delta_*(1)
\end{equation}
and
\begin{equation}\label{d*-}
\delta_*(a) =\delta^{-1}(\delta(1)a\delta(1)), \qquad  a\in {\cal A}
\end{equation}
where $\delta^{-1}\!:\delta({\cal A})\to P\cal A$ is the inverse
mapping to\ \ $\delta\!:P\cal A\to \delta(\cal A)$.
\end{Tm}

\begin{I}\label{coeff}
Let $A\subset L(H)$ be a $^*$-subalgebra containing the identity 1
of $L(H)$ and $V\in L(H)$. We call $A$ the {\em coefficient
algebra of the $C^*$-algebra $C^* (A,V)$} generated by $A$ and $V$
if  $A$  and   $V$ satisfy the following three conditions
\begin{gather}\label{c1}
Va= VAV^*V, \qquad a \in A;\\[6pt]
\label{c2}
VaV^* \in A, \qquad a \in A\\
\intertext{and}
\label{c3}
V^* aV\in A, \qquad a\in A.
\end{gather}

It was shown in \cite{Leb-Odz} and \cite{Bakht-Leb} that instead of
condition \eqref{c1} one can use equivalently the condition
\begin{equation}
\label{c4}
 V^*V\in Z(A)
\end{equation}
where $Z(A)$ is the center of $A$ or the condition
\begin{equation}
\label{c8}
 V^*VaV^* bV = aV^*bV, \qquad a,b \in A.
\end{equation}
It is worth mentioning that conditions \eqref{c1}, \eqref{c2} and
\eqref{c3} imply that $V$ is a partial isometry and the mapping $A\ni
a \mapsto VaV^*$ is an endomorphism \cite[Proposition
2.2]{Leb-Odz}. Thus (recalling \eqref{b,,2}) we see that a
$C^*$-algebra $A\subset L(H)$ containing the identity of $L(H)$ is the
coefficient algebra for $C^*(A,V)$ iff the mapping $V\,\cdot\,V^* \! :
A\to A$ is an endomorphism and the mapping $V^*\,\cdot\,V \! : A\to A$ is
a transfer operator for $V\,\cdot\,V^*$.
\end{I}

\begin{I}\label{coef-d}
Let $\delta$ be an endomorphism of an (abstract) unital $C^*$-algebra
$\cal A$. We say that the pair $(\cal A,\delta)$ is
{\em finely representable\/} if there exists a triple
$(H,\pi ,U)$ consisting of a Hilbert space~$H$, faithful
non-degenerate representation $\pi\!:{\cal A}\to L(H)$ and a linear
continuous operator \hbox{$U\!:H\to H$} such that for every $a\in\cal A$ the
following conditions are satisfied
\begin{gather}\label{b,,4}
\pi(\delta(a)) =U\pi(a)U^*,\qquad U^*\pi(a)U \in \pi(\cal A)\\
\intertext{and}
\label{sdvig}
U\pi (a) = \pi(\delta(a))U, \qquad a \in \cal A.
\end{gather}
That is $\pi ({\cal A})$ is the coefficient algebra for $C^*(\pi
({\cal A}), U)$ under the  fixed  endomorphism $U \,\cdot\,U^*$.
In this case we also say that $\cal A$ is a {\em coefficient algebra
associated with $\delta$.}

By the foregoing discussion (see \ref{coeff}) instead of condition
\eqref{sdvig} one can use equivalently the condition
\begin{equation}
\label{c11} U^*U\in Z(\pi ({\cal A}))
\end{equation}
or the condition
\begin{equation} \label{c9}
 U^*U\pi (a)U^* \pi (b)U = \pi (a)U^*\pi (b)U, \ \ a,b \in \cal A.
\end{equation}

In particular it is clear that the finely representable pair $(\cal
A,\delta)$ can also be defined as a pair such that there exists a
triple $(H,\pi ,U)$ where $\pi\!:{\cal A}\to L(H)$ is a faithful
non-degenerate representation, $U\in L(H)$ and the mapping $U\,\cdot\,
U^*$ coincides with the endomorphism $\delta$ on $\pi (\cal A)$ while
the mapping $U^*\,\cdot\,U $ is a transfer operator for $\delta$.

Since $\delta$ is an endomorphism it follows that $\delta (1)$ is a
projection and so \eqref{b,,4} implies that $UU^*$ is a projection, so
$U$ is a partial isometry.
\end{I}

\begin{Tm}
\label{b..1} {\upshape \cite[Theorem  3.1]{Bakht-Leb} }\ \
A pair\/ $(\cal A,\delta)$ is finely representable iff there
exists a complete transfer operator $\delta_*$ for $({\cal A},\delta)$.
\end{Tm}

\begin{Dn}\label{crossed}\ \
Let $(\cal A,\delta)$ be a finely representable pair. The {\em
crossed product of\/} $({\cal A},\delta)$, which we denote by ${\cal
A}\times_{\delta} {\mathbb Z}$, or simply by ${\cal A}\times {\mathbb
Z}$ when $\delta$ is understood, is the universal unital $C^*$-algebra
generated by a copy of $\cal A$ and a partial isometry $U$ subject to
relations
\begin{equation}
\label{,b,,4}
\delta(a) = UaU^*,\quad\ \delta_* (a)=U^*aU ,\qquad a\in \cal A,
\end{equation}
where $\delta_*$ is the complete transfer operator for $(\cal A,
\delta)$ (this $\delta_*$ does exist by Theorem \ref{b..1} and is
unique by Theorem \ref{complete}). The algebra $\cal A$ will be
called the {\em coefficient algebra\/} for $\cal A\times_\delta
\mathbb Z$.
\end{Dn}

{\em Remark}. The reason why we use $\mathbb Z$ in the notation of
the crossed product but not $\mathbb N$ (as in a number of
sources) will be uncovered in what follows (in the next section).

\medbreak
Theorems \ref{complete} and \ref{b..1} imply the non-degeneracy of
Definition \ref{crossed} (there exists a non-zero representation for
${\cal A}\times_{\delta}{\mathbb Z}$ and $\cal A$ is a
$C^*$-subalgebra in ${\cal A}\times_{\delta} {\mathbb Z}$).
The further investigation of the structure of $\cal A\times_\delta
\mathbb Z$ is presented if the next section.

\section{Faithful and regular representations of\hfil\null\ \hbox{the crossed product}\hfill}
\label{faithful}

This section is devoted to the description of the internal structure
of ${\cal A}\times_{\delta} {\mathbb Z}$. Among the main technical
instruments here are the results of \cite{Leb-Odz} and
\cite{Bakht-Leb}.

\begin{I}\label{2.1}
Let $\hat{a}$ and $\hat{U}$ be the canonical
images of\ \ $a\in \cal A$ and $U$ in  ${\cal A}\times {\mathbb
Z}$, respectively. Note that we can identify $\hat{a}$ with $a$
which justifies the usage of notation $\delta (\hat{a})$ and
$\delta_* (\hat{a})$.

By the definition of the crossed product and \ref{coef-d} we have that
$\hat{\cal A}$ is a coefficient algebra of ${\cal A}\times_{\delta}
{\mathbb Z} = C^* (\hat{\cal A}, \hat{U})$ (see \ref{coeff}).
\end{I}

\begin{Pn}\label{B0'}{\upshape \cite[Proposition 2.3]{Leb-Odz} }\ \
Let $A$ be the coefficient algebra of $C^* (A,V)$ (see \ref{coeff}).
Then the
vector space $B_0$ consisting of finite sums
\begin{equation}\label{suma'}
x= V^{*N}a_{-N}+ \dots + V^*a_{-1}+ a_0 +
a_1V +\dots + a_NV^N ,
\end{equation}
where $a_k\in A$ and $N\in{\mathbb N}\cup\{0\}$, is a
dense $^*$-subalgebra of  $C^* (A,V)$.
\end{Pn}

Since any $C^*$-algebra can be faithfully represented as a
$C^*$-subalgebra of operators acting in some Hilbert space
Proposition \ref{B0'} implies  the next

\begin{Pn}\label{B0}
The vector space\/ $C_0$ consisting of finite sums
\begin{equation}\label{suma}
x= \hat{U}^{*N}\hat{a}_{-N}+ \dots +
\hat{U}^*\hat{a}_{-1}+ \hat{a}_0 + \hat{a}_1\hat{U} +\dots +
\hat{a}_N\hat{U}^N ,
\end{equation}
where $a_k\in\cal A$ and $N\in{\mathbb N}\cup\{0\}$, is a
dense $^*$-subalgebra of ${\cal A}\times_{\delta}{\mathbb Z}$.
\end{Pn}

\begin{I}\label{*'}\,
We  say that the algebra $C^*(A,V)$ mentioned
in \ref{coeff} possesses the {\em property\/} (*) if for
any $x\in B_0$ given by \eqref{suma'} the inequality
\begin{equation} \label{gwiazdka}
\norm{a_0}\le \norm{x} \qquad\qquad (^*)
\end{equation}
holds.
\end{I}

\begin{I}
\label{2.6} Observe that  ${\cal A}\times_{\delta}
{\mathbb Z} = C^* (\hat{\cal A}, \hat{U})$ possesses  property
(*). Indeed, take any  faithful non-degenerate representation
$\pi\!:{\cal A}\to  L(H)$ and a partial isometry
$U\!:H\to H$ mentioned in \ref{coef-d}.
Consider the space ${\cal H} =l^2 ({\mathbb Z}, H)$ and the
representation  $\nu \! : C^* (\hat{\cal A}, \hat{U}) \to L({\cal
H})$ given by the formulae
\begin{gather*}
(\nu (\hat{a})\xi )_n = \pi (a) (\xi_n), \qquad\text{where}\quad
a\in {\cal A}, \quad l^2 ({\mathbb Z}, H) \ni \xi = \{ \xi_n  \}_{n\in
{\mathbb Z}}\,;\\[6pt]
(\nu (\hat{U})\xi )_n = U (\xi_{n-1}),\qquad
(\nu (\hat{U}^*)\xi )_n = U^* (\xi_{n+1}).
\end{gather*}
\end{I}
Routine verification shows that $\nu (\hat{\cal A})$ and $\nu
(\hat{U})$ satisfy all the conditions mentioned in \ref{coef-d}
(for $\pi ({\cal A})$ and $U$).

Now take any $x\in  C^* (\hat{\cal A}, \hat{U})$ given by
\eqref{suma} and for a given $\varepsilon > 0$ chose  a vector
$\eta \in H$ such that
\begin{equation}\label{e}
\Vert \eta \Vert =1 \quad \text{and}\quad \Vert \pi (a_0)
\eta \Vert > \Vert \pi (a_0)  \Vert - \varepsilon .
\end{equation}
Set $\xi \in l^2 ({\mathbb Z}, H)$ by $\xi_n = \delta_{0n}\eta $,
where $\delta_{ij}$ is the Kronecker symbol. We have that $\Vert
\xi \Vert = 1$ and the explicit form of $\nu (x) \xi$  and
\eqref{e} imply
\[
\Vert \nu (x) \xi \Vert^2 \ge \Vert \pi (a_0) \eta \Vert^2 >
(\Vert \pi (a_0)  \Vert - \varepsilon )^2
\]
which by the arbitrariness of $\varepsilon$ proves the desired
inequality
\[
\Vert x  \Vert \ge \Vert \pi (a_0)  \Vert = \Vert  \hat a_0 \Vert.
\]

The foregoing observation provides  us a possibility to exploit
all the main results of \cite{Leb-Odz} to uncover the structure of
the crossed product and we start to do this.

The next theorem shows  the crucial value of property
(*) in the crossed product --- this property is a criterium for
a representation of the crossed product to be faithful.

\begin{Tm}\label{iso-crossed*}\ \
Let $\cal A$ be a $C^*$-algebra with an
identity\/ $1$ and $\delta\!: {\cal A} \to {\cal A}$ be an
endomorphism such that there exists a complete transfer operator
$\delta_*$ for $({\cal A}, \delta )$ and let $(H,\pi ,U)$ be a
triple mentioned in \ref{coef-d}. Then the map
\begin{equation}
\begin{split}
&\label{is1} \Phi (\hat{a}) = \pi (a), \qquad a\in {\cal A},\\[6pt]
&\qquad\quad \Phi (\hat{U}) = U
\end{split}
\end{equation}
gives rise to the isomorphism between the algebras $ {\cal
A}\times_{\delta} {\mathbb Z}= C^* (\hat{\cal A},
\hat{U})$ and $ C^* (\pi ({\cal A}), U) $ iff the algebra $C^*
(\pi ({\cal A}), U)$ possesses property\/ $({^*})$.
\end{Tm}
{\bfseries Proof.}\ \
Follows from \ref{2.6} and \cite[Theorem 2.13]{Leb-Odz}.
\qed

\begin{Rk}\label{faith}
Note that the representation $\nu \! : C^*
(\hat{\cal A}, \hat{U}) \to L({\cal H})$ described in \ref{2.6}
possesses property (*) so it is a faithful representation of
${\cal A}\times_{\delta} {\mathbb Z}$.
\end{Rk}
\begin{Rk}\label{cros}
It is worth mentioning that the value of property (*) in the theory of
crossed products of $C^*-$algebras by discrete groups (semigroups) of
automorphisms (endomorphisms) has been observed by a number of authors
along with the proofs of the results of Theorem \ref{iso-crossed*}
type. The importance of property (*) for the first time (probably) was
clarified by O'Donovan \cite{O'Donov} in connection with the
description of $C^*-$algebras generated by weighted shifts. The most
general result establishing the crucial role of this property in the
theory of crossed products of $C^*-$algebras by discrete groups of
{\em automorphisms} was obtained in \cite{Leb1} (see also
\cite{AntLeb}, Chapters 2, 3 for complete proofs and various
applications) for an {\em arbitrary} $C^*-$algebra and {\em amenable}
discrete group. The relation of the corresponding property to the
faithful representations of crossed products by {\em endomorphisms}
generated by {\em isometries} was investigated in \cite{BKR,ALNR}. The
role of the properties of (*) type in the theory of Fell bundles
algebras was studied in \cite{Exel-amenab}. The properties of this
sort proved to be of great value not only in pure $C^*-$theory but
also in various applications such as, for example, the construction of
symbolic calculus and developing the solvability theory of functional
differential equations (see \cite{AnLebBel1, AnLebBel2}).

We shall also exploit this property heavily in the subsequent part
of the paper.
\end{Rk}

\begin{I}
\label{regular} {\bfseries Regular representation of the crossed
product}.\ \
The representation mentioned in Remark \ref{faith} being faithful is
not defined explicitly (in terms of ${\cal A}$, $\delta$, $\delta_*$), in
fact we have only established its existence. Now we shall present a
faithful representation of ${\cal A}\times_{\delta} {\mathbb Z}$ that
will be written out explicitly in terms of ${\cal A}$, $\delta$,
$\delta_*$. Keeping in mind the standard regular representations for
the known various versions of crossed products it is reasonable to
call it the {\em regular representation of\/} $ {\cal A}\times_\delta
{\mathbb Z}$. In fact the construction of this
representation has been obtained in the proof of Theorem 3.1 in
\cite{Bakht-Leb}.

First we construct the desired Hilbert $H$ space by means of the
elements of the initial algebra~$\cal A$ in the following way. Let
$\langle\,\cdot\,,\,\cdot\,\rangle$ be a certain non-negative inner
product on $\cal A$ (differing from a common inner product only in
such a way that for certain non-zero elements $v\in\cal A$ the
expression $\langle v,v\rangle$ may be equal to zero). For example
this inner product may have the form $\langle v,u\rangle =f(u^*v)$
where $f$ is some positive linear functional on $\cal A$. If one
factorizes $\cal A$ by all the elements $v$ such that $\langle
v,v\rangle =0$ then he obtains a linear space with a strictly positive
inner product. We shall call the completion of this space with respect
to the norm $\norm v =\sqrt{\langle v,v \rangle}$ the {\em Hilbert
space generated by the inner product} $\langle\,\cdot\,,\,
\cdot\,\rangle$.

We shall build a desired triple $(H,U,\pi)$ by means of the triple
$(\cal A,\delta,\delta_*)$. Let $F$ be the set of all positive linear
functionals on $\cal A$. The space $H$ will be constructed as the
completion of the direct sum $\bigoplus_{f\in F}H^f$ of some Hilbert
spaces $H^f$. Every $H^f$ will in turn be the completion of the direct
sum of Hilbert spaces $\bigoplus_{n\in\mathbb Z} H^f_n$. These $H^f_n$
are generated by non-negative inner products
$\langle\,\cdot\,,\,\cdot\,\rangle_n$ on the initial algebra $\cal A$
that are given by the following formulae
\begin{align}
\label{b,,7}
\langle v,u\rangle_0 &=f(u^*v);\\[3pt]
\label{b,,8} \langle v,u\rangle_n &=f\bigl(\delta_*^n(u^*v)\bigr),
\qquad n\ge 0;\\[3pt]
\label{b,,9} \langle v,u\rangle_n
&=f\bigl(u^*\delta^{|n|}(1)v\bigr), \qquad n\le 0.
\end{align}
Note the next   properties of these inner products
(\cite{Bakht-Leb}, Lemma 3.2): for any \/ $v,u\in\cal A$ the
following equalities are true
\begin{align}
\label{b,,10} \langle \delta(v),u\rangle_{n+1} &=\langle
v,\delta_*(u)\rangle_n, \qquad n\ge 0;\\[3pt]
\label{b,,11} \langle\delta^{|n|}(1)v,u\rangle_{n+1} &=\langle
v,\delta^{|n|}(1)u \rangle_n,\qquad n<0.
\end{align}

Now let us define the operators  $U$ and $U^*$ on the space  $H$
constructed. These operators  leave invariant all the subspaces
$H^f\subset H$. The action of  $U$ and $U^*$ on every $H^f$ is the
same and its scheme is presented in the first line of the next
diagram.

\noindent
\begin{alignat*}{14}
&\dots\ \,&&
\xrlarrow[\,]{\delta^3(1)\,\cdot\,}{\delta^3(1)\,\cdot\,}&& \
H^f_{-2}\ &&
\xrlarrow[\,]{\delta^2(1)\,\cdot\,}{\delta^2(1)\,\cdot\,}&& \
H^f_{-1}\ &&
\xrlarrow[\,]{\delta(1)\,\cdot\,}{\delta(1)\,\cdot\,}&& \ H^f_{0}\
&& \xrlarrow[\,]{\delta}{\delta_*}&& \ H^f_{1}\ &&
\xrlarrow[\,]{\delta}{\delta_*}&& \ H^f_{2}\ &&
\xrlarrow[\,]{\delta}{\delta_*}& \ \,\dots&&&
\qquad \xrlarrow[\ ]{\textstyle U}{\textstyle\,U^*}\\[6pt]
&\dots &&&&\,\mspace{0mu}\delta^2(a) &&&&\ \delta(a) &&&&\ \;a
&&&&\ \;a &&&&\ \;a &&& \dots &&&\qquad\ \:\pi(a)
\end{alignat*}

\medbreak\noindent Formally this action is defined in the
following way. Consider any {\em finite} sum
\[
h =\bigoplus_n h_n \in H^f, \qquad h_n\in H^f_n .
\]
Set
\[
Uh =\bigoplus_n(Uh)_n \quad \text{and} \quad U^*h =\bigoplus_n(U^*h)_n
\]
where
\begin{align}
\label{b,,12} (Uh)_n &=
\begin{cases}
\delta(h_{n-1}),&\ \ \text{if}\ \ n>0,\\
\delta^{|n|+1}(1)h_{n-1},& \ \ \text{if}\ \ n\le 0,
\end{cases}\\[6pt]
(U^*h)_n &=
\begin{cases}
\delta_*(h_{n+1}), & \ \ \text{if}\ \ n\ge 0,\\
\delta^{|n|}(1)h_{n+1},&\ \  \text{if}\ \ n<0.
\end{cases}
\label{b,,13}
\end{align}
Equalities \eqref{b,,10} and \eqref{b,,11}
guarantee that the operators $U$ and $U^*$ are
well defined  (i.\,e.\ they preserve factorization and completion
by means of which the spaces $H^f_n$ were built from the algebra
$\cal A$) and $U$ and $U^*$ are mutually adjoint.

Now let us define the representation  $\pi\!:{\cal A}\to L(H)$.
For any $a\in\cal A$ the operator  $\pi(a)\!:H\to H$ will leave
invariant all the subspaces $H^f\subset H$ and also all the
subspaces  $H^f_n \subset H^f$. If  $h_n\in H^f_n$ then we set
\begin{equation}
\label{b,,14} \pi(a)h_n =
\begin{cases}
ah_n, & \ \ n\ge 0,\\
\delta^{|n|}(a)h_n, &\ \ n\le 0.
\end{cases}
\end{equation}
The scheme of the action of the operator  $\pi(a)$ is presented in
the second line of the diagram given above.

In the process of the proof of Theorem 3.1 in \cite{Bakht-Leb}
there was verified that the triple $(H, \pi , U)$ described above
satisfies all the conditions of \ref{coef-d} and property (*) here
can be proved in the same way as in  \ref{2.6}.
\end{I}

The results of \cite{Leb-Odz} provide us an opportunity to say
more on the structure and properties of elements in $ {\cal
A}\times_{\delta} {\mathbb Z}$ and henceforth we
present these properties.

\begin{I}\label{a_k}\ \
Since $\hat{U} \,\cdot\,\hat{U}^* \! : \hat{\cal A}\to \hat{\cal A}$
is an endomorphism it follows that for any $k\in {\mathbb N}$ the map
${\hat{U}}^k \,\cdot\,{\hat{U}}^{*k} \! : \hat{\cal A}\to \hat{\cal
A}$ is an endomorphism as well and therefore ${\hat{U}}^k
{\hat{U}}^{*k}$ is a projection and thus $\hat U^k$, \,$k\in {\mathbb N}$
is a partial isometry. So $\hat{U}^k=\hat{U}^k\hat{U}^{*k}\hat{U}^k$
and $\hat{U}^{*k}=\hat{U}^{*k}\hat{U}^k\hat{U}^{*k}$. Therefore one
can always choose coefficients $\hat{a}_k$ and $\hat{a}_{-k}$ of
\eqref{suma} in such a way that
\begin{equation}
\label{wybor}
\hat{a}_k\hat{U}^k\hat{U}^{*k}=\hat{a}_k\qquad\textrm{and }\qquad
\hat{U}^k\hat{U}^{*k}\hat{a}_{-k}=\hat{a}_{-k}.
\end{equation}
The forthcoming proposition shows that the assumption that
$\hat{a}_k,\hat{a}_{-k}\in \hat{\cal A}$, $k=1,\,\ldots,\,N$ satisfy
\eqref{wybor} guarantees their uniqueness in the expansion
\eqref{suma} and more over by means of these coefficients any element
of $ {\cal A}\times_{\delta} {\mathbb Z}$ can be determined in a
unique way. This is the subject of Theorem \ref{uniqueNk}.
\end{I}

\begin{Pn} \label{gwiazdka2p}
If the   coefficients of $x$ in \eqref{suma} satisfy
\eqref{wybor} then
\begin{equation}
\label{gwiazdka2} \norm{\hat{a}_k} =\norm{a_k}\le \norm x,
\qquad \norm{\hat{a}_{-k}}=\norm{a_{-k}}\le \norm{x}
\end{equation}
for $k\in\{0,1,\ldots,N\}$. And in particular these
coefficients are uniquelly defined.
\end{Pn}
{\bfseries Proof}. Follows from  \cite[Proposition 2.6]{Leb-Odz} since
${\cal A}\times_{\delta} {\mathbb Z}$  possesses property (*).
\qed

\begin{I}\label{N_k}\ \
Proposition \ref{gwiazdka2p}  means that  one
can define the linear and continuous maps ${\cal N}_k \!: C_0 \to
{\cal A}$ and ${\cal N}_{-k} \!: C_0 \to {\cal A}$, \,
$k\in\mathbb N$
\begin{equation}
\begin{split}\label{N-k}
&\,{\cal N}_k(x)=a_k\in{\cal A} U^kU^{*k}\subset\cal A,\\[6pt]
&{\cal N}_{-k}(x)=a_{-k}\in U^kU^{*k}{\cal A}\subset{\cal A}.
\end{split}
\end{equation}
By continuity these mappings can be expanded onto the whole of
${\cal A}\times_{\delta} {\mathbb Z}$ thus defining the
{\bfseries 'coefficients'} of an arbitrary element $x\in  {\cal
A}\times_{\delta} {\mathbb Z}$.

Theorem \ref{3a.N}  presented  below shows that  the norm of an
element $x\in C_0$ can be calculated only in terms of the elements
of ${\cal A}$ (0-degree coefficients of the powers of $xx^*$).
\end{I}

\begin{Tm}\label{3a.N}
For any element $x\in C_0$ of the form \eqref{suma} we have
\begin{equation}\label{be3.131}
\Vert x \Vert = \lim_{k\to\infty}
\sqrt[\leftroot{-2}\uproot{1}\scriptstyle 4k]{
\left\Vert {\cal N}_0 \left[ (xx^*)^{2k}\right]\right\Vert }
\end{equation}
where ${\cal N}_0$ is the mapping defined by \eqref{N-k}.
\end{Tm}
{\bfseries Proof}. Follows from   \cite[Theorem 2.11]{Leb-Odz}.
\qed

\medbreak
We finish the section with the statement showing that the mapping
${\cal N}_0$ is an exact conditional expectation from $  {\cal
A}\times_{\delta} {\mathbb Z}$ onto ${\cal A}$ and
any element $x\in {\cal
A}\times_{\delta} {\mathbb Z}$ can be 'restored' by its
 coefficients ${\cal N}_k (x)$ and ${\cal N}_{-k} (x)$.

\begin{Tm}\label{uniqueNk}\ \
Let\/  $x\in {\cal A}\times_{\delta} {\mathbb Z}$.
Then the following conditions are equivalent:

\smallskip
\quad\llap{$(i)$}\ \ $x=0;$

\smallskip
\quad\llap{$(ii)$}\ \ ${\cal N}_k (x)=0$, \,$k\in\mathbb Z;$

\smallskip
\quad\llap{$(iii)$}\ \ ${\cal N}_0 (x^*x)=0$.
\end{Tm}
{\bfseries Proof}. Follows from \cite[Theorem 2.15]{Leb-Odz}.
\qed

\section{Various crossed  products}
\label{v-crossed product}

In this section we list a number of successful constructions of
crossed products associated with this or that sort of $C^*$-algebras
and endomorphism and discuss their interrelations with the crossed
product introduced.

\subsection{Monomorphisms with hereditary range}

We start with the crossed-product under special hypotheses which have
been considered in \cite{Murphy}.

Let  $\cal A$ be a $C^*$-algebra with an identity $1$ and
$\delta\!:\cal A\to \cal A$ be a  monomorphism  with hereditary
range. We shall denote by ${\cal U}({\cal A}, \delta )$ the universal
unital $C^*$-algebra generated by ${\cal A}$ and an isometry $T$
subject to the relation
\begin{equation}\label{T}
\delta (a) = TaT^*, \qquad a\in {\cal A}.
\end{equation}
This ${\cal U}({\cal A}, \delta)$ has been proposed in \cite{Murphy}
as the definition for the crossed-product of $A$ by $\delta$.

Since $\delta (\cal A )$ is a hereditary subalgebra of $\cal A$ it
follows \cite[Proposition 4.1]{exel} that $\delta ({\cal A} )= \delta
(1){\cal A} \delta (1)$. So we find ourselves under the conditions
of \,3) of Theorem \ref{complete} with $P=1$. Therefore for the
endomorphism $\delta$ under consideration there is a unique complete
transfer operator $\delta_*$ given by \eqref{d*-}. Thus we can take
the crossed product $ {\cal A}\times_{\delta} {\mathbb Z}$. The
definition \ref{crossed} of the crossed product $ {\cal A}
\times_{\delta} {\mathbb Z}$, the universal property of ${\cal U}
({\cal A}, \delta)$ and \eqref{T} implies that $ {\cal A}
\times_{\delta} {\mathbb Z}$ is a representation of ${\cal U}({\cal
A}, \delta)$ and since ${\cal A}$ is embedded in $ {\cal A}
\times_{\delta} {\mathbb Z}$ it is embedded in ${\cal U}({\cal A},
\delta)$ as well. We shall therefore view ${\cal A}$ as a subalgebra
of ${\cal U}({\cal A}, \delta)$.

\begin{Pn}\label{monom}
Let  $\cal A$ be a $C^*$-algebra with an identity
$1$ and $\delta\!:\cal A\to \cal A$ be a  monomorphism  with
hereditary range. Then the mapping
\[
\varphi \!: {\cal U}({\cal A}, \delta) \to {\cal A}\times_{\delta}
{\mathbb Z}
\]
such that $\varphi (T) = \hat{U}$ and $\varphi (a) = \hat{a}$ for all
$a\in {\cal A}$, where $\hat{a}$ and $\hat{U}$ are the canonical
images of\ \ $a\in \cal A$ and $U$ in ${\cal A}\times_{\delta}
{\mathbb Z}$, respectively, is a $^*$-isomorphism.
\end{Pn}
{\bfseries Proof}.
The argument preceding the proposition  imply that
$\varphi $ is a $^*$-epimorphism and to finish the proof it is
enough to observe that
\[
T^*aT =\delta_* (a), \qquad a\in \cal A .
\]
But this follows from \eqref{T}, \,\eqref{d*-}, equality $T^*T=1$
and the relations
\[
T^*aT =T^*TT^* aTT^*T =T^*\delta (1)a\delta (1)T =T^*\delta
\delta_* (a)T =T^*T\delta_* (a)T^*T =\delta_* (a). \qed
\]

\subsection{Partial crossed product}
\label{PCP}

The notion of the partial crossed product of a $C^*$-algebra by a
group $\mathbb Z$ of partial automorphisms was introduced by R.\ Exel
\cite{Exel1}, it was generalized by K.\ McClanahan
\cite{McClanahan} up to the crossed product of a $C^*$-algebra by
partial actions of discrete groups. Faithful representations of
these crossed products were described in \cite{top-free}.

Within the framework of the present article it is natural  to
confine ourselves to the group~$\mathbb Z$. Let us recall the notion of
the partial crossed product \cite{Exel1}.

\begin{I}\label{p-c}
A {\em partial automorphism} of a $C^*$-algebra $ A$ is a triple
$\Theta =(\theta , I,J)$ where $I$ and $J$ are closed two sided ideals
in $ A$ and $\theta \!: I \to J$ is a $^*$-isomorphism.

Given a partial automorphism one can consider, for each integer $n$,
the set $D_n$ --- the domain of $\theta^{-n}$ (which is equal to the
image of $\theta^{n}$). By convention we put $D_0 =A$ and $\theta^0$
is the identity automorphism of $A$.

For each integer $n, \,D_n$ is an ideal in $A$ \cite[Proposition
3.2]{Exel1}.

For the sake of convenience of the presentation we give hereafter the
definition of the partial crossed product which differs from the
original Exel's definition but in fact is equivalent to it (cf.
\cite[Sections 3 and 5]{Exel1}).

The {\em partial crossed product} for $\Theta =(\theta , I,J)$ is the
universal enveloping $C^*$-algebra $C^* (A, \Theta )$ generated by
finite sums
\begin{equation}\label{e.p-cr.1}
a_{-n}V^{*n}+ \dots+ a_{-1}V^{*} + a_0 + a_1V +\dots  a_nV^n, \qquad
n\in \mathbb N
\end{equation}
where $a_i\in D_i$, \,$i\in \mathbb Z$ and $V$ is a partial isometry
such that

(i) \,the initial space of $V$ is $I \cal H$ (meaning closed linear
span) and the final space of $V$ is $J\cal H$ (here $\cal H$ is the
space where $A$ acts), and

(ii) \,$Va V^* = \theta (a), \ \ a\in I$.

\begin{Rk}\label{R1}
Condition (i) implies in particular that both the projections $V^*V$
and $VV^*$ belong to the commutant of $A$. In addition conditions (i)
and (ii) imply that $V^*a V = \theta^{-1} (a)$, \,$a\in J$.
\end{Rk}
\end{I}

The next proposition describes the conditions on a coefficient
algebra under which the crossed product introduced in this article
is a partial crossed product.

\begin{Pn}\label{partial aut}
Let $\delta_*$ be a complete transfer operator for $({\cal A},
\delta)$ and $\delta (1)\in Z({\cal A})$, then

\indent $1)$ \,the triple $\Theta = (\delta , \delta_* (1){\cal A},
\delta (1) {\cal A})$ is a partial automorphism, and

\indent $2)$ \,$\delta_* \!: {\cal A} \to {\cal A}$ is an
endomorphism, and

\indent $3)$ \,the mapping
\[
\phi \!: C^* ({\cal A}, \Theta ) \to {\cal A}\times_{\delta}\mathbb Z
\]
such that $\varphi (V) = \hat{U}$ and $\varphi (a) = \hat{a}$ for
all $a\in {\cal A}$, where $\hat{a}$ and $\hat{U}$ are the
canonical images of\ \ $a\in \cal A$ and $U$ in  ${\cal
A}\times_{\delta} {\mathbb Z}$, respectively, is a
$^*$-isomorphism.
\end{Pn}
{\bfseries Proof.} 1) and 2) are  proved in \cite[Proposition
2.9]{Bakht-Leb}.

By applying equality \eqref{b,,3} and the condition $\delta (1)\in
Z({\cal A})$ we obtain for any $a,b \in \cal A$
\[
\delta (\delta_* (a)b) = \delta \delta_* (a)\delta (b)= \delta
(1)a\delta (1)\delta (b) = a \delta (b).
\]
This means that $\delta$ is a transfer operator for $\delta_*$.
Therefore we have for every $n\in \mathbb N$ \,$\delta ^n (1)\in
Z({\cal A})$ and $\delta_* ^n (1)\in Z({\cal A})$, and it is easy to
see that $D_n = {\cal A} \delta^n (1)$, and $D_{-n} = {\cal A}
\delta_*^n (1)$ (here $D_n$ are the ideals mentioned in \ref{p-c}).

Now the result follows from the definition of $ C^* ({\cal A}, \Theta
) $ along with the definition of \,${\cal A}\times_{\delta}\mathbb Z$ and
Proposition \ref{B0} (here one should note that since $\delta^n (1)=
\hat{U}^n \hat{U}^{*n}\in Z(\hat{\cal A})$ we have $\hat{U}^{*n}
\hat{a} = \hat{U}^{*n}\hat{U}^n \hat{U}^{*n} \hat{a}= \hat{U}^{*n}
\hat{a}\hat{U}^n \hat{U}^{*n}= \delta_*^n (\hat{a})\hat{U}^{*n}
\hat{U}^n \hat{U}^{*n} = \delta_*^n (\hat{a})\delta_*^n
(1)\hat{U}^{*n}$, and $\hat{a}\hat{U}^n= \hat{a}\hat{U}^n
\hat{U}^{*n}\hat{U}^n = \hat{a}\delta^n (1)\hat{U}^n$, and therefore
the form of \eqref{suma} coincides with the form of \eqref{e.p-cr.1}).
\qed

\medbreak
Proposition \ref{partial aut} is 'almost invertible'. That is by a
'slight' (natural) extension of a partial automorphism one can
always obtain a coefficient algebra satisfying the conditions of
Proposition~\ref{partial aut}.

Indeed. Let $ A$ be a unital algebra and $\Theta =(\theta , I,J)$
be those mentioned in \ref{p-c}. Set $A_1 = \{ A, VV^*, V^*V \}$
to be the (enveloping) $C^*$-algebra generated by $A$, $VV^*$ and
$V^*V$ ($V$ is the partial isometry mentioned in the definition
of the partial crossed product). The sets $VV^*A$ and $V^*VA$ are
the ideals in $A_1$ and the isomorphism $\theta \!: I \to J$ extends
up to the isomorphism
\[
\tilde{\theta} \!: V^*VA \to VV^*A.
\]
Evidently the mapping
\[
\delta \!: A_1 \to A_1, \qquad \delta (\,\cdot\,) = \tilde{\theta}(V^*V
\,\cdot \,)
\]
is an endomorphism of $A_1$ and the mapping
\[
\delta_* \!: A_1 \to A_1, \qquad \delta_* (\,\cdot\,) =
\tilde{\theta}^{-1}(VV^* \,\cdot \,)
\]
is a complete transfer operator for $(A_1,\delta)$.

Since $\delta(1)= VV^* \in Z(A_1) $ it follows that $(A_1, \delta,
\delta_*)$ satisfy all the conditions of Proposition~\ref{partial aut}.

Thus by a slight abuse of language one can say that the partial
crossed product is the crossed product introduced in this article
under the additional condition $\delta (1) \in Z({\cal A})$.

\subsection{Cuntz-Krieger algebras}
\label{Cuntz}

Throughout this subsection we shall let $A$ be an $n\times n$ matrix
with $A(i,j)\in\{0,1\}$ for all $i$ and~$j$, and such that no row and
no column of $A$ is identically zero. The Cuntz-Krieger algebra ${\cal
O}_A$ (see \cite{CK}) is a $C^*$-algebra generated by partial
isometries $S_i$, \,$i\in \overline{1,n}$ that act on a Hilbert space
in such a way that their support projections $Q_i =S^*_iS_i$ and their
range projections $P_i = S_iS^*_i$ satisfy the relations
\begin{equation}
\label{C-K} P_iP_j =0, \qquad i\neq j;\qquad\qquad
Q_i = \sum_{r=1}^n A(i,r)P_r , \qquad  i,j \in \overline{1,n}\,.
\end{equation}
The algebras of this sort arise naturally as the objects associated
with topological Markov chains and serve as a source of inspiration
for numerous investigation.

We shall show that ${\cal O}_A$ can be naturally considered as a
certain  crossed product of the type introduced.

\medbreak
The sum of the range projections $P_i$ is a unit in  ${\cal O}_A$.
If $\mu = (i_1,\dots,i_k)$ is a multiindex with $i_j \in
\overline{1,n}$ we denote by $\vert \mu \vert$ the length $k$ of
$\mu$ and write $S_\emptyset =1$, \,$S_\mu = S_{i_1}S_{i_2}\dotsm
S_{i_k}$ ($\emptyset$ is also considered as a multiindex). It is
shown in \cite{CK} that all $S_\mu$ are partial isometries.

We recall in this connection that the product of two partial
isometries is not necessarily a partial isometry. The criterium
for the  product of partial isometries to be a partial isometry is
given in the next

\begin{Pn}\label{prod-partis}{\upshape\cite[Lemma 2]{Hal-Wal}}\ \
Let $S$ and
$T$ be partial isometries. Then $ST$ is a partial isometry iff
$S^*S$ commutes with  $TT^*$.
\end{Pn}
The symbols $P_\mu$, $Q_\mu$ will stand for the range and support
projections of $S_\mu$, respectively. The foregoing observation means
that for any two given multiindexes $\mu$ and $\nu$ the projections
$P_\mu $ and $ Q_\nu$ commute.

In addition by  Lemma 2.1 of \cite{CK}, for $\vert\mu \vert =\vert
\nu\vert $ we have
\begin{equation}\label{SmuSnu}
S_\mu^*S_\nu \neq 0\ \Rightarrow\ \mu=\nu \quad \text{and}
\quad S_\mu^*S_\nu =Q_\mu =Q_{i_k},\quad (\mu =(i_1,\dotsc,i_k))
\end{equation}
which implies
\begin{equation}
\label{a(Pij)} P_\mu P_\nu = \delta _{\mu ,\nu }P_\mu
\end{equation}
and since $\sum_i P_i =1$ it also follows that
\begin{equation}
\label{(P=1)} \sum_{\vert \mu \vert =k}P_\mu =1.
\end{equation}

Let ${\cal F}_A$ be the $C^*$-algebra generated by all the
elements of the form $S_\mu P_iS^*_\nu$, where
$\vert\mu\vert=\vert\nu\vert =k$, \ $k=0,\,1,\,\dotsc$; \
$i\in \overline{1,n}$. Clearly
\begin{equation}\label{subset}
S_i {\cal F}_A S_i^* \subset {\cal F}_A, \qquad
i\in \overline{1,n}
\end{equation}
and it is also shown in \cite[Lemma 2.2]{CK}  that
\begin{equation}
\label{subset1} S_i^* {\cal F}_A S_i \subset {\cal F}_A, \qquad
i\in \overline{1,n}.
\end{equation}

Moreover  \cite[Proposition 2.8]{CK}  tells  that any element $X$
of the $^*$-algebra $\cal O_A$ generated by $S_i$, \,$i\in
\overline{1,n}$ \,can be written as a finite sum
\begin{equation}
\label{CK*} X= \sum_{\vert \nu\vert \ge 1} X_\nu S^*_\nu +X_0 +
\sum_{\vert \mu\vert \geq 1}S_\mu X_\mu
\end{equation}
where $X_\nu,X_0,X_\mu \in {\cal F}_A$ and for $X$ in
\eqref{CK*} the next inequality is true
\begin{equation}\label{CK**}
\Vert X \Vert \ge \Vert X_0 \Vert.
\end{equation}
For every projection $P_i$,  \,$i\in \overline{1,n}$ we set
\begin{equation}
\label{a(P)} \alpha (P_i)= \sum_{j=1}^n S_jP_iS_j^* = \sum_{j=1}^n
P_{j,i}.
\end{equation}
Clearly $\alpha (P_i)\in {\cal F}_A$. Moreover
\eqref{a(P)}, \eqref{a(Pij)} and \eqref{(P=1)} imply
\begin{equation}
\label{a(Pi)} \alpha (P_i)\alpha (P_j)= 0, \quad i\neq j
\end{equation}
and
\begin{equation}
\label{a(P=1)} \sum_{i}\alpha (P_i) =1.
\end{equation}
We extend by linearity  formula \eqref{a(P)} to
\begin{equation}\label{a(Qi)}
\alpha (Q_i ) = \sum_{r=1}^n A(i,r)\alpha ( P_r ).
\end{equation}
Set
\begin{equation}\label{Q}
Q = \sum_{i=1}^n Q_i
\end{equation}
and
\begin{equation}\label{a(Q)}
\alpha (Q)=\sum_{i=1}^n \alpha (Q_i ).
\end{equation}
Obviously $\alpha (Q_i ), Q, \alpha (Q) \in {\cal F}_A$ and since
no column of the matrix $A$ is zero it follows from \eqref{C-K}
along with \eqref{Q} and from \eqref{a(P=1)}, \eqref{a(Qi)} and
\eqref{a(Q)} that
\begin{equation}
 \label{a(Q)>1}
Q \ge 1 , \qquad \alpha (Q) \ge 1 \,.
\end{equation}
Set
\begin{equation}\label{S}
S =\alpha (Q)^{-\frac{1}{2}}\sum_{i=1}^n S_i.
\end{equation}
From this and  \eqref{SmuSnu} we  get
\begin{equation}
\label{Si} S_i = P_i\alpha (Q)^{\frac{1}{2}}S.
\end{equation}

Observe now that $S$ is an isometry.
Indeed, by definition of $Q$ and $\alpha (Q)$, \eqref{Q},
\eqref{a(Q)} and \eqref{a(Qi)} we have that
\begin{equation}\label{aQ1}
\alpha (Q) = \sum_{r=1}^n \gamma_r\alpha (P_r),
\end{equation}
where $\gamma_r \neq 0$, \,$r=\overline{1,n}$ and
\begin{equation}\label{Q1}
 Q = \sum_{r=1}^n \gamma_r P_r ,
\end{equation}
and therefore
\begin{equation}\label{aQ-1}
\alpha (Q)^{-1} = \sum_{r=1}^n \gamma_r^{-1}\alpha (P_r),
\end{equation}
and
\begin{equation}\label{Q-1}
Q^{-1} = \sum_{r=1}^n \gamma_r^{-1} P_r .
\end{equation}
Now we have by \eqref{S} and \eqref{SmuSnu}
\[
S^*S = \sum_{i,j} S_i^*\alpha (Q)^{-1}S_j = \sum_{j} S_j^*\alpha
(Q)^{-1}S_j
\]
and substituting \eqref{aQ-1} into the latter formula and
recalling \eqref{a(P)}, \eqref{Q-1} and \eqref{SmuSnu} we obtain
\begin{gather*}
S^*S =\sum_{j} S_j^*\sum_r \gamma_r^{-1}\alpha (P_r)S_j
=\sum_{j}\sum_r\gamma_r^{-1}\sum_iS_j^*S_iP_rS_i^*S_j\\[6pt]
=\sum_{j}\sum_r\gamma_r^{-1}Q_jP_r=
\sum_{j}Q_j\sum_r\gamma_r^{-1}P_r =Q Q^{-1}=1 .
\end{gather*}
Thus $S$ is an isometry.

In view of \eqref{S},  \eqref{subset} and \eqref{subset1} we
conclude that
\begin{equation}\label{F-coefficient}
S{\cal F}_A S^* \subset {\cal F}_A \quad
\text{and }\quad S^*{\cal F}_A S \subset {\cal F}_A .
\end{equation}
Since $S$ is an isometry it follows that the mapping $\delta \!:
{\cal F}_A \to {\cal F}_A$ defined by
\begin{equation}\label{dS}
\delta (\,\cdot\,) = S (\,\cdot\,)S^*
\end{equation}
is an endomorphism of ${\cal F}_A$.

Observe that \eqref{F-coefficient} and \eqref{dS} mean that ${\cal
F}_A$ is a coefficient algebra for the $C^*$-algebra $C^* ({\cal
F}_A,S)$ and
\begin{equation}\label{dS*}
\delta_* (\,\cdot\,) = S^* (\,\cdot\,)S .
\end{equation}

Now we are ready to establish the desired crossed product
structure of the Cuntz-Krieger $C^*$-algebra ${\cal O}_A$.

\begin{Pn}
\label{CK-cr} For every $n\times n$ matrix $A$ with no zero rows
and columns  we have
\[
{\cal O}_A = C^* ({\cal F}_A,S)\cong {\cal F}_A\times_{\delta}{\mathbb Z}
\]
where the isomorphism
\[
\varphi \!: C^* ({\cal F}_A,S)\to {\cal F}_A\times_{\delta} {\mathbb Z}
\]
is such that $\varphi (S) = \hat{U}$ and $\varphi (a) = \hat{a}$ for
all $a\in {\cal F}_A$, where $\hat{a}$ and $\hat{U}$ are the canonical
images of\ \ $a\in {\cal F}_A$ and $U$ in ${\cal F}_A\times_{\delta}
{\mathbb Z}$, respectively and $\delta$ and $\delta_*$ are defined by
\eqref{dS} and \eqref{dS*}.
\end{Pn}
{\bfseries Proof.}
The equality ${\cal O}_A = C^* ({\cal F}_A,S)$ follows from \eqref{S},
\eqref{Si} along with \eqref{CK*} and observation that $P_i$, \,$i\in
\overline{1,n}$ and $\alpha (Q)$ belong to ${\cal F}_A$.

Note now that \eqref{CK**} is nothing else than property (*) for
$C^* ({\cal F}_A,S)$. Therefore the desired result follows from
Theorem \ref{iso-crossed*}.
\qed

\subsection{Paschke's type crossed product}

While analyzing the simplicity of the Cuntz algebra (considered in
\cite{cuntz}) W.\,L.\ Paschke has found a certain condition on the
action of an endomorphism generated by an isometry $S$ under which the
$C^*$-algebra $C^*(A,S)$ generated by the initial $C^*$-algebra $A$
and $S$ is on the one hand isomorphic to a certain crossed product and
on the other hand is simple. His result is stated as follows.

\begin{Tm}\label{pas0}{\upshape \cite[Theorem 1]{Paschke}}
Let $ A$ be a strongly amenable unital $C^*$-algebra acting on a
Hilbert space $H$. Suppose that $S$ is a nonunitary isometry (i.\,e.\
$S^*S=1\neq SS^*$) in $L(H)$ such that

\smallskip
\quad\llap{$(i)$}\ \,$SAS^* \subset A, \ S^*AS \subset A$; and

\smallskip
\quad\llap{$(ii)$}\ \,if a proper (two-sided) ideal $J$
satisfies the condition $SJS^* \subseteq J$ then $J =\{0\}$.

\smallskip\noindent
Then $C^*(A,S)$ is simple.
\end{Tm}
We shall observe now that Paschke's result can be easily
generalized up to the situation considered in this article.
First we note the next

\begin{La}\label{pas2}
Let $A$ be the coefficient algebra of \,$C^* (A,V)$ \,(see
\ref{coeff}). Suppose that $A$ is strongly amenable, $VV^*\neq 1$
and  the only proper (two-sided) ideal $J$ for which $VJV^*
\subseteq J$ is the zero ideal. Then $C^* (A,V)$ possesses
property \emph{(*)}.
\end{La}
{\bfseries Proof.}\ \
The proof can be obtained by the word by word repetition of the
corresponding parts of the proofs of \cite[Lemmas 2 and 3]{Paschke}
using when necessary instead of the Paschke's condition $S^*S=1$ the
coefficient algebra condition $V^*V\in Z(A)$. So we omit it.
\qed

\medbreak
As a corollary of this lemma we obtain the following generalization of
Theorem \ref{pas0}.

\begin{Tm}\label{pas-gener}
Let $A$ be the coefficient algebra of \,$C^* (A,V)$ \,(see
\ref{coeff}). Suppose that $A$ is strongly amenable, $VV^*\neq 1$ and
if a proper (two-sided) ideal $J$ satisfies the condition $VJV^*
\subseteq J$ then it is zero. Then

\smallskip
$1)$\ \,$C^*(A,V) \cong A\times_{\delta} {\mathbb Z}$
\ where the isomorphism
\[
\varphi \!: C^*(A,V) \to A\times_{\delta} {\mathbb Z}
\]
is such that $\varphi (V) = \hat{U}$ and $\varphi (a) = \hat{a}$ for
all $a\in A$, where $\hat{a}$ and $\hat{U}$ are the canonical images
of $a\in A$ and $U$ in $A\times_{\delta} {\mathbb Z}$, respectively
and $\delta (\,\cdot\,) =V(\,\cdot\,)V^*$ and $\delta_* (\,\cdot\,)
= V^* (\,\cdot\,)V$; and

\smallskip
$2)$\ \,$C^*(A,V)$ is simple.
\end{Tm}
{\bfseries Proof.} \,1) By Lemma \ref{pas2} $C^*(A,V)$ possesses
property (*). Thus the result follows from Theorem \ref{iso-crossed*}.

\medbreak
2) \,Let $I$ be any proper ideal in $C^*(A,V)$. Then $J=I\cap A$ is a
proper ideal in $A$ (since if $J=A$ then $I = C^*(A,V)$). Clearly
$VJV^* \subseteq J$ and therefore
\begin{equation}\label{J}
I\cap A=J=\{ 0\}.
\end{equation}
Consider the canonical mapping  $\pi \!: C^*(A,V)\to C^*(A,V)/I$.
It is easy to see that  $\pi (C^*(A,V))= C^* (\pi (A), \pi (V))$,
and for every $a\in A$ we have
\begin{gather*}
\pi (V^*aV) = \pi (V)^* \pi (a) \pi
(V),\quad \text{and}\quad \pi (VaV^*) = \pi (V) \pi (a) \pi (V)^*, \quad
\text{and}\\[6pt]
\pi (V^*V)\in Z(\pi (A)).
\end{gather*}
The equality \eqref{J} implies
\[
\pi (A) \cong A/A\cap I\cong A
\]
and therefore $\pi (C^*(A,V))= C^* (\pi (A), \pi (V))$ satisfies
all the assumptions of Lemma \ref{pas2} with $A$, $V$ substituted for
$\pi (A)$, $\pi (V)$. Thus by the already proved first part of the
theorem we have that
\[
\pi (C^*(A,V))\cong  A\times_{\delta} {\mathbb Z}\cong C^*(A,V).
\]
And it follows that $I=\{ 0\}$.
\qed

\subsection{Exel's crossed product}

Recently R.\ Exel proposed in \cite{exel} a new definition for the
crossed product of a unital $C^*$-algebra by an endomorphism
$\alpha $ and a certain transfer operator $\cal L$. He also proved
in \cite{exel} that this new construction generalizes many of the
previously known constructions among which are  the above
mentioned monomorphisms with hereditary range, Paschke's crossed
product and Cuntz-Krieger algebras (in fact we believe that the
main inspiring motivation  for   Exel was to generalize the
Cuntz-Krieger construction, and he perfectly succeeded).

In this subsection we analyze the interrelations between Exel's
crossed product and ours. We show that on the one hand the crossed
product introduced in this article is a special  case of Exel's
and on the other hand in the most natural situations (when all the
powers of $\cal L$ are generated by partial isometries) Exel's
crossed product is of the type introduced here but with {\em
different} algebra, {\em different} endomorphism
and {\em different} transfer operator. In fact in the
general  situation the crossed product 'philosophy' is even more
peculiar and subtle and this will be the theme of the subsequent
sections.

Throughout this subsection $A$ will be a unital $C^*$-algebra and
$\alpha \!:A\to A$ will be a \hbox{*-endo}\-morphism and $\cal L$ will
be a certain transfer operator for the pair $(A,\alpha)$, that is
${\cal L} \!: A \to A$ is a continuous positive map satisfying
condition \eqref{b,,2} with $\delta $ substituted for $\alpha $ and
$\delta_*$ for $\cal L$.

\begin{I}\label{Ex}
The (Exel's) crossed product $A\times_{\alpha , \cal L}{\mathbb N}$
is the universal $C^*$-algebra generated by a copy of $A$ and
element $S$ subject to the relations

\smallskip
\quad\llap{(i)}\ \ $Sa = \alpha (a)S, \quad a\in A$,

\smallskip
\quad\llap{(ii)}\ \ $S^* aS = {\cal L} (a), \ \ a\in A$,

\smallskip
\quad\llap{(iii)}\ \ if $(a,k)\in \overline{A\alpha (A) A}\times
\overline{ASS^*A}$ is such that
\begin{equation}\label{redun}
abS = kbS
\end{equation}
for all $b\in A$ then $a=k$ (here $\overline{A\alpha (A) A}$ is the
closed linear space generated by $A\alpha (A) A$ (the two-sided ideal
in $A$ generated by $\alpha (A)$) and $\overline{ASS^*A}$ is the
closed linear space generated by $ASS^*A$).

Any pair $(a,k)$ mentioned in (iii) is called a  {\em redundancy}.
\end{I}

\begin{I}\label{Tau}
Let ${\cal T} (A, \alpha , {\cal L})$ be the universal algebra
satisfying (i) and (ii). Exel showed \cite[3.5]{exel} that for {\em
any} endomorphism $\alpha$ and a transfer operator $\cal L$ the
algebra ${\cal T} (A, \alpha , {\cal L})$ is non-degenerate (i.\,e.\
the canonical map
\[
A \to {\cal T} (A, \alpha , {\cal L})
\]
is injective).

Therefore one can think of the crossed product $A\times_{\alpha ,
\cal L}{\mathbb N}$ as the quotient of ${\cal T} (A, \alpha ,
{\cal L})$ by the closed two-sided ideal $\cal I$ generated by the
set of differences $a-k$, for all redundancies $(a,k)$.

In contrast to the situation with ${\cal T} (A, \alpha , {\cal
L})$  not for all $\alpha$ and $\cal L$ there is a natural
inclusion of $A$ in $A\times_{\alpha , \cal L}{\mathbb N}$. The
conditions that ensure the canonical map $A\to A\times_{\alpha ,
\cal L}{\mathbb N}$ to be injective were found by N.\ Brownlowe and
I.\ Raeburn  \cite{Brow-Rae}. They are formulated in terms of
Cuntz-Pimzner algebras.
\end{I}

We start with the observation that the crossed product introduced in
this paper is a special case of the Exel's crossed product. This is
stated in Theorem \ref{My=Ex}.

The argument here goes absolutely in the same way as the
corresponding argument in Section 4 of \cite{exel} and we give the
necessary  proofs simply for the sake of completeness of
presentation.

Evidently ${\cal A}\times_{\delta} {\mathbb Z}$ is not 'bigger' than
${\cal T} ({\cal A}, \delta , \delta_*)$ which is stated in the next

\begin{Pn}\label{3.4}
Let the pair\/ $(\cal A,\delta)$ be finely representable and\/
$\delta_*$ be the corresponding transfer operator.
Then there exists a unique *-epimorphism
\[
\psi\!:  {\cal T} ({\cal A}, \delta , \delta_*) \to {\cal
A}\times_{\delta} {\mathbb Z}
\]
such that $\psi(S) = \hat U$, and $\psi(a)=\hat a$, for all $a\in
{\cal A}$.
\end{Pn}
{\bfseries Proof.}\ \
Since $\delta (\,\cdot\,) = U(\,\cdot\,)U^*$ and $U^*U \in Z
({\cal A})$ it follows  that, for all $a\in \cal A$,
\[
U a = UU^*U a  = U a U^*U = \delta (a)U.
\]
In other words, relations \ref{Ex} (i), (ii) hold for $\hat {\cal A}$
and $\hat U$ within ${\cal A}\times_{\delta} {\mathbb Z}$
and hence the conclusion follows from the universal property of ${\cal
T} ({\cal A}, \delta , \delta_*)$.
\qed

\medbreak
To establish the coincidence between  ${\cal A}\times_{\delta}
{\mathbb Z}$ (the crossed product introduced in this
article) and ${\cal A}\times_{\delta , \delta_*} {\mathbb N}$ (the
Exel's crossed product)  we observe first the following fact.

\begin{Pn}\label{SIsAnIsometry}\ \
Let $({\cal A}, \delta ,\delta_*)$ satisfy the hypothesis of
Proposition \ref{3.4}. Then:

\smallskip
\quad\llap{$(i)$}\ \,The canonical element $S\in {\cal T} ({\cal A},
\delta , \delta_*)$ is a partial isometry, and hence also its image
$\dot S\in {\cal A}\times_{\delta , \delta_*}{\mathbb N}$.

\smallskip
\quad\llap{$(ii)$}\ \,For every $a\in {\cal A}$ one has that
$(\delta (a), S a S^*)$ is a redundancy.
\end{Pn}
{\bfseries Proof.}\ \ (i)\ \ By \eqref{P} we have that $\delta_*(1)$
is a projection thus $S$ is a partial isometry.

(ii)\ \ For any  $b\in {\cal A}$ we have
\[
SaS^*bS = Sa\delta_* (b) = \delta (a)\delta \delta_*(b)S =
\delta (a)\delta (1)b\delta (1)S = \delta(a)bS,
\]
where we have used the completeness of $\delta_*$ and the fact that $
\delta(1) S = S 1 = S$. It remains to notice that $SaS^* = \delta
(a)SS^*\in \overline {{\cal A} S S^* {\cal A}}$.
\qed

\begin{Cy}\label{3.2}
Under the hypotheses of Proposition \ref{SIsAnIsometry}
there exists a unique *-epi\-morphism
\[
\phi\!: {\cal A}\times_{\delta}{\mathbb Z} \to
{\cal A}\times_{\delta , \delta_*}{\mathbb N}
\]
such that $\phi(\hat{U})=\dot S$, and $\phi(\hat{a})=\dot a$, for all
$a\in \cal A$, where $\dot S$ and $\dot a$ are the canonical images of
$S$ and $a$ in ${\cal A}\times_{\delta , \delta_*}{\mathbb N}$,
respectively and $\hat{a}$ and $\hat{U}$ are the canonical images of
$a\in \cal A$ and $U$ in ${\cal A}\times_{\delta} {\mathbb Z}$ (here
${\cal A}\times_{\delta , \delta_*}{\mathbb N}$ is the Exel's crossed
product and ${\cal A}\times_{\delta} {\mathbb Z}$ is the crossed
product introduced in this paper).
\end{Cy}
{\bfseries Proof.}
By Proposition \ref{SIsAnIsometry} (ii) we have that $\delta (\dot a)
= \dot S \dot a \dot S^*$. Hence the conclusion follows from the
universal property of ${\cal A}\times_{\delta} {\mathbb Z}$.
\qed

\medbreak
Thus ${\cal A}\times_{\delta,\delta_*}{\mathbb N}$ is not 'bigger'
than ${\cal A}\times_{\delta} {\mathbb Z}$. The next result shows the
desired coincidence between ${\cal A}\times_{\delta} {\mathbb Z}$ and
${\cal A}\times_{\delta , \delta_*} {\mathbb N}$.
\begin{Tm}\label{My=Ex}
Let\/ $(\cal A,\delta)$ be a finely representable pair with the
complete transfer operator\/ $\delta_*$. Then the map $\phi$ of
Corollary \ref{3.2} is a *-isomorphism between ${\cal
A}\times_{\delta} {\mathbb Z}$ and ${\cal A}\times_{\delta , \delta_*}
{\mathbb N}$.
\end{Tm}
{\bfseries Proof. }\ \
Observe first that the map $\psi$ of Proposition \ref{3.4} vanishes on
the ideal $\cal I$ mentioned in \ref{Tau}. Indeed. Let $(a,k)\in
\overline{{\cal A} \delta ({\cal A}) {\cal A}}\times \overline{{\cal
A}SS^*{\cal A}}$ be a redundancy. Therefore for all $b\in {\cal A}$
one has $abS = kbS$. By applying of $\psi$ to both sides one obtains
\[
\hat {a}\hat{ b}\hat{ U} = \psi(k) \hat {b} \hat {U}.
\]
Since $\delta (1) = UU^*$ we have for all $b,c\in {\cal A}$ that
\[
\hat a\hat b\widehat{\delta (1)}\hat c =  \hat a\hat b\hat U\hat U^*\hat c =
\psi(k)\hat b\hat U\hat U^*\hat c = \psi(k)\hat b\widehat{\delta(1)}\hat c.
\]
It follows that $\hat {a}\hat{ x} = \psi(k)\hat {x}$ for all $ x\in
\overline {{\cal A}\delta (1){\cal A}}$. Since $k\in \overline{{\cal
A}SS^*{\cal A}}$ we have that $\psi(k)\in \overline{\hat{{\cal
A}}\hat{ U}\hat{U}^*\hat{{\cal A}}} = \overline{\hat{{\cal A}}
\widehat{\delta (1)} \hat{{\cal A}}}$. Finally by the completeness of
$\delta_*$ we have
\[
a \in  \overline{{\cal A} \delta ({\cal A}) {\cal A}} =
\overline {{\cal A} \delta (1) {\cal A}\delta (1) {\cal A}} =
\overline {{\cal A}\delta (1) {\cal A}}.
\]
Thus $\hat{a}=\psi(k)$ and it follows that $\psi(a-k)=0$ and hence
that $\psi$ vanishes on $\cal I$ as claimed. By passage to the
quotient we get a map
\[
\widetilde\psi \!: {\cal A}\times_{\delta , \delta_*} {\mathbb N}\to
{\cal A}\times_{\delta} {\mathbb Z},
\]
which is the inverse of the map $\phi$ of Corollary \ref{3.2}.
\qed

\medbreak
Thus we have established that the crossed product introduced in the
present article is a special case of the Exel's one. Now we shall move
in the opposite direction and show that in 'the most popular'
situation when all the operators $S^n$, \,$n=1,\,2,\,\dotsc$ are
partial isometries (which is equivalent to the condition that all the
elements ${\cal L}^n (1)$, \,$n=1,\,2,\,\dotsc$ are projections) the
Exel's crossed product is of type introduced here but with {\em
different} algebra, {\em different} endomorphism and {\em different}
transfer operator.

We start with the result showing that under the mentioned hypothesis
even ${\cal T} ( A, \alpha , {\cal L})$ possesses the structure of a
certain crossed product.
\begin{Tm}\label{Tau-cr}
Let\/ $A$ be a unital $C^*$-algebra, $\alpha \!: A \to A$ be an
endomorphism, and $\cal L$ be a transfer operator such that all the
elements ${\cal L}^n (1)$, \,$n=1,\,2,\,\dotsc$ are projections. Let
$\cal A$ be the $C^*$-algebra generated by $A$ and $S^kS^{*k}$,
\,$k=1,\,2,\,\dotsc$ where $S$ is the universal operator satisfying
relations $(i)$ and $(ii)$ of \ref{Ex}. Then the map
\[
\nu \!: {\cal A}\times_{\delta} {\mathbb Z} \to
{\cal T}(A,\alpha,{\cal L})
\]
such that $\nu (\hat{U}) = (S)$, and $\nu (\hat{a}) = a$, for all
$a\in A$, and $\delta \!: {\cal A} \to \cal A$ is given by $\delta
(\,\cdot\,) = S (\,\cdot\,)S^*$, and $\delta_* \!: {\cal A} \to \cal
A$ is given by $\delta_* (\,\cdot\,) = S^* (\,\cdot\,)S$ establishes
a $^*$-isomorphism.
\end{Tm}
{\bfseries Proof.}\ \
As it has been already mentioned the condition that all the elements
${\cal L}^n (1)$, \,$n=1,\,2,\,\dotsc$ are projections is equivalent
to the condition that all the operators $S^k$, \,$k=1,\,2,\,\dotsc$
are partial isometries. In view of Proposition \ref{prod-partis} this
condition implies that all the operators $S^{*k}S^k$,\ \ $S^jS^{*j}$,\
\ $k,j =1,\,2,\,\dotsc$ commute with each other and therefore $S^*S$
belongs to the commutant of $\cal A$. Observe now that for any $a\in
A$ and any $k=1,\,2,\,\dotsc$ we have
\[
SaS^* = \alpha (a)SS^*\in {\cal A}, \ \ \text{and}\ \ S(S^kS^{*k})S^* =
S^{k+1}S^{*k+1} \in {\cal A}
\]
and for any $c,d \in {\cal A}$ we have
\begin{equation}\label{S0}
ScdS^* = SS^*ScdS^*= ScS^*SdS^*
\end{equation}
where we have used the mentioned fact that $S^*S$ belongs to the
commutant of $\cal A$.

The foregoing observation implies
\begin{equation}\label{S1}
S{\cal A}S^* \subset {\cal A} \quad \text{and}\quad S^*S \in Z({\cal A}).
\end{equation}
Note also that \eqref{S0} implies that $S(\,\cdot\,)S^*$ is an
endomorphism of $\cal A$.

A routine computation shows also that
\begin{equation}\label{S2}
S^*{\cal A}S\subset \cal A .
\end{equation}

In view of the definition of ${\cal T} ( A, \alpha , {\cal L})$ and
$\cal A$ we have
\begin{equation}\label{S3}
{\cal T} ( A, \alpha , {\cal L}) = C^* ({\cal A}, S).
\end{equation}
And \eqref{S1} and \eqref{S2} mean that $\cal A$ is the coefficient
algebra for $C^* ({\cal A}, S)$ with $\delta (\,\cdot\,)= S(\,\cdot\,)
S^*$ and $\delta_* (\,\cdot\,) = S^* (\,\cdot\, )S$. Now the
universal property of ${\cal A}\times_{\delta} {\mathbb Z}$ implies
that the mapping
\[
{\nu} \!: {\cal A}\times_{\delta} {\mathbb Z} \to
{\cal T} ( A, \alpha , {\cal L})
\]
where ${\nu} ({\hat a})=a$, \,$a\in A$ and ${\nu} (\hat{U})=S$ is a
$^*$-epimorphism.

To prove that this mapping is in fact a $^*$-isomorphism it is enough
to show that the $C^*$-alge\-bra $C^* ({\cal A}, S)$ in \eqref{S3}
possesses property (*) (in this case one can apply Theorem
\ref{iso-crossed*}). So let us verify the latter property. Let $A$ and
$S$ be the universal algebra and element that generate ${\cal T} ( A,
\alpha , {\cal L}) = C^* ({\cal A}, S)$. They satisfy relations
\eqref{S1} and \eqref{S2}, and without loss of generality we can
assume that $C^* ({\cal A}, S)$ is a $C^*$-subalgebra of $ L(H)$ for
some Hilbert space $H$ and that the identity of $A$ is the identity of
$L(H)$. Consider the space ${\cal H} =l^2 ({\mathbb Z}, H)$ and the
representation $\mu\!: {\cal T} ( A, \alpha , {\cal L}) \to
L({\cal H})$ given by the formulae
\begin{gather*}
(\mu ({a})\xi )_n =  a (\xi_n), \qquad  a\in { A}, \quad l^2
({\mathbb Z}, H) \ni \xi = \{ \xi_n  \}_{n\in {\mathbb Z}}\,;\\[6pt]
(\mu (S)\xi )_n = S (\xi_{n-1}),\qquad (\mu {S}^*)\xi )_n =
S^* (\xi_{n+1}).
\end{gather*}
It is easy to see that $\mu ({\cal A})\cong {\cal A} $, and that $\mu
({\cal A})$ and $\mu (S)$ satisfy the same relations \eqref{S1} and
\eqref{S2}, and that $\mu (S)$ generates the same mappings $\delta$
and $\delta_*$ on $\mu ({\cal A})$. Note now that the algebra $C^*
(\mu ({\cal A}), \mu (S))$ possesses property (*) (recall the argument
of \ref{2.6}). But this means that the algebra $C^* ({\cal A},
S)={\cal T} ( A, \alpha , {\cal L})$ possesses this property as well.
The proof is finished.
\qed

\medbreak
Now we establish the desired isomorphism between the Exel's
crossed product and the crossed product introduced in this
article.

\begin{Tm}\label{Ex!}
Let the hypothesis of Theorem \ref{Tau-cr} be satisfied. Let
$\dot{\cal A}$ be the $C^*$-algebra generated by $\dot A$ and ${\dot
S}^k{\dot S}^{*k}$, \ $k=1,\,2,\,\dotsc$ where $\dot S$ and $\dot A$
are the canonical images of $S$ and $A$ in $A\times_{\alpha , \cal
L}{\mathbb N}$, respectively. Then the the map
\[
{\gamma} \!: {\dot{\cal A}}\times_{\delta} {\mathbb Z} \to
A\times_{\alpha , \cal L}{\mathbb N}
\]
such that $\gamma (\hat{U}) = \dot S$, and $\gamma (\hat{a}) = \dot
a$, for all $a\in A$, and $\delta \!: \dot{\cal A} \to \dot{\cal A}$
is given by $\delta (\,\cdot\, ) = \dot S (\,\cdot\, ){\dot S}^*$, and
$\delta_* \!: \dot{\cal A} \to \dot {\cal A}$ is given by $\delta_*
(\,\cdot\, ) = {\dot S}^* (\,\cdot\, )\dot S$ establishes a
$^*$-isomorphism.
\end{Tm}
{\bfseries Proof.}
In view of \eqref{S1}, \eqref{S2} and \eqref{S3} we have that
\begin{equation}\label{e-c}
A\times_{\alpha , \cal L}{\mathbb N} = C^* (\dot{\cal A}, \dot S)
\end{equation}
and $\dot {\cal A}$ is the coefficient algebra for $C^* (\dot {\cal
A}, \dot S)$ with $\delta (\,\cdot\, )= \dot S (\,\cdot\, ){\dot S}^*$
and $\delta_* (\,\cdot\, ) = {\dot S}^* (\,\cdot\, )\dot S$. Now the
universal property of $\dot {\cal A}\times_{\delta} {\mathbb Z}$
implies that the mapping
\[
{\gamma} \!: {\dot {\cal A}}\times_{\delta} {\mathbb Z} \to
A\times_{\alpha , \cal L}{\mathbb N}
\]
where ${\gamma} ({\hat a})=\dot a, \ a\in A$ and ${\gamma}
(\hat{U})=\dot S$ is a $^*$-epimorphism.

To prove that this mapping is in fact a $^*$-isomorphism it is enough
to show that the \hbox{$C^*$-}algebra $C^* (\dot{\cal A}, \dot S)$ in
\eqref{e-c} possesses property~(*).

Let $\dot A$ and $\dot S$ be the universal algebra and element that
generate $A\times_{\alpha , \cal L}{\mathbb N} = C^* (\dot{\cal A},
\dot S)$. They satisfy relations \eqref{S1} and \eqref{S2} (for
$\dot{\cal A}$ and $\dot S$), and without loss of generality we can
assume that $C^* (\dot{\cal A}, \dot S)$ is a $C^*$-subalgebra of $
L(H)$ for some Hilbert space $H$ and that the identity of $\dot A$ is
the identity of $L(H)$. Consider the space ${\cal H} =l^2 ({\mathbb
Z}, H)$ and the representation $\nu \! : {\cal T} ( A, \alpha , {\cal
L}) \to L({\cal H})$ given by the formulae
\begin{gather*}
(\nu ({a})\xi )_n =  {\dot a} (\xi_n), \qquad  a\in { A}, \quad l^2
({\mathbb Z}, H) \ni \xi = \{ \xi_n  \}_{n\in {\mathbb Z}}\,;\\[6pt]
(\nu (S)\xi )_n = {\dot S} (\xi_{n-1}),\qquad (\nu ({S}^*)\xi)_n =
{\dot S}^* (\xi_{n+1}).
\end{gather*}
It is easy to see that $\nu ({\cal A})\cong \dot {\cal A} $ (here
${\cal A}$ is the $C^*$-algebra mentioned in Theorem \ref{Tau-cr}),
and that $\nu ({\cal A})$ and $\nu (S)$ satisfy the same relations
\eqref{S1} and \eqref{S2} (for $\dot{\cal A}$ and $\dot S$), and that
$\nu (S)$ generates the same mappings $\delta$ and $\delta_*$ on $\nu
({\cal A})$. Moreover since $\dot A$ and $\dot S$ satisfy relations
(i), (ii), and (iii) of \ref{Ex} then by the construction of $\nu$ we
have that $\nu (A)$ and $\nu (S)$ satisfy these relations as well.
Therefore we can consider $\nu$ as a representation of
$A\times_{\alpha , \cal L}{\mathbb N}= C^* (\dot{\cal A}, \dot S)$.
But $C^* (\nu ({\cal A}), \nu (S))$ possesses property $(^*)$ (by the
argument of \ref{2.6}). Thus $C^* (\dot{\cal A}, \dot S)$ possesses
this property as well.
\qed

\subsection{Kwasniewski's  crossed product}
\label{Kwa}

This subsection is devoted to the description of the crossed
product structure developed recently by B.\,K.\ Kwasniewski
\cite{kwa} and the discussion of the interrelation between this
structure and the crossed product introduced in the present
paper. To give the motivation of the Kwasniewski's crossed product
we start with two simple examples.

\begin{Ee}\label{1}
Consider the Hilbert space $H=L^2({\mathbb R})$. Let
$A\subset L(H)$ be the $C^*$-algebra of operators of
multiplication by continuous bounded functions on $\mathbb R$ that
are constant on ${\mathbb R}_- = \{ x: x\le 0 \}$. Set the unitary
operator $U\in L(H)$ by the formula
\begin{equation}\label{U1}
(Uf)(x)= f (x-1), \qquad f (\,\cdot\, )\in H .
\end{equation}
Routine verification shows that the mapping
\begin{equation}\label{U2}
A\ni a \mapsto UaU^*
\end{equation}
is the  endomorphism of $A$ of the form
\begin{equation}\label{U3}
UaU^* (x) = a(x-1), \qquad a(\,\cdot\, ) \in A ,
\end{equation}
and
\begin{equation}\label{U4}
U^*aU (x) = a(x+1), \qquad a(\,\cdot\, ) \in A .
\end{equation}
Clearly the mapping $A\ni a \mapsto U^*aU$ {\em does not} preserve $A$.

Let $C^* (A,U)$ be the $C^*$-algebra generated by $A$ and $U$. It is
easy to show that
\begin{equation}
\label{U5} C^* (A,U) = C^* ({\cal A},U)
\end{equation}
where ${\cal A}\subset L(H)$ is the $C^*$-algebra of operators of
multiplication by continuous bounded functions on $\mathbb R$ that
have limits at $-\infty$.

In addition we have that
\begin{equation}
\label{U6} U{\cal A}U^*\subset {\cal A} \qquad \text{and}\qquad
U^*{\cal A}U\subset {\cal A}
\end{equation}
and the corresponding actions $\delta (\,\cdot\, )= U(\,\cdot\, )U^*$
and $\delta_* (\,\cdot\, )= U^* (\,\cdot\, )U$ on ${\cal A}$ are given
by formulae \eqref{U3} and \eqref{U4}.

Thus ${\cal A}$ is a coefficient algebra for $C^* ({\cal A},U)$.

Moreover one can easily check that $C^* ({\cal A},U)$ possesses
property $(*)$. Therefore by  Theorem \ref{iso-crossed*} we
conclude that
\begin{equation}\label{U7}
C^* (A,U) = C^* ({\cal A},U)\cong {\cal A}\times_{\delta}{\mathbb Z}.
\end{equation}
We would like to emphasize that the situation considered {\em does
not} satisfy the conditions under which the {\em Exel's} crossed
product was defined since we started with the algebra $A$ and the
operator $U$ such that the mapping $A\ni a \mapsto U^*aU$ \eqref{U4}
is {\em not} a transfer operator (it does not preserve $A$). But after
extending $A$ up to ${\cal A}$ we have obtained the coefficient
algebra and thus found ourselves under the main assumptions of the
crossed product construction of the present article.
\end{Ee}

In fact by a slight modification of the example considered it is
easy to 'worsen' the situation even further.
\begin{Ee}\label{22}
Let $H$ and $U$ be the same as in the previous example, and
${A}\subset L(H)$ be the $C^*$-algebra of operators of multiplication
by continuous bounded functions on $\mathbb R$ that are constant on
${\mathbb R}_- = \{ x: x\le 0 \}$ and are constant for $x\ge \pi$.
Then we have that
\begin{equation}
\label{U8} C^* ({A}, U) = C^* ({{\cal A}}, U),
\end{equation}
where ${{\cal A}}$ is the $C^*$-algebra of operators of multiplication
by continuous bounded functions on $\mathbb R$ that have limits at
$\pm \infty$.

As in the previous example here ${\cal A}$ is the coefficient
algebra for $C^* ({{\cal A}}, U)$ and  we have  that
\[
C^* (A,U) = C^* ({\cal A},U)\cong {\cal A}\times_{\delta}{\mathbb Z}
\]
where $\delta$ and $\delta_*$ are given by the same formulae. But in
contrast to the previous example here even the mapping $A\ni a \mapsto
UaU^*$ {\em does not} preserve $A$ (thus this mapping is {\em not\/} an
endomorphism).
\end{Ee}

These examples show that neither a transfer operator nor even a
certain endomorphism are among the starting objects that lead to the
coefficient algebras and crossed products. In fact in both these
examples the principal moment was a certain procedure of {\em
extension} of the initial algebra $A$ up to a coefficient algebra
${\cal A}$ (this procedure is generated by the mappings $A\ni a
\mapsto UaU^*$ and $A\ni a \mapsto U^*aU$). After the implementation
of this procedure and obtention of the coefficient algebra ${\cal A}$
the final step --- the construction of the crossed product goes
smoothly (in accordance with the scheme introduced in the present
article).

The general procedures of extension of initial $C^*$-algebras up to
coefficient algebras were given in \cite{Leb-Odz}. The maximal ideal
spaces of the arising in this way {\em commutative} coefficient
algebras were described in \cite{Kwa-Leb}. By developing the technique
of \cite{Leb-Odz} and \cite{Kwa-Leb} along with a number of new ideas
B.\,K.\ Kwasniewski has described the extension procedure as for the
initial {\em commutative} $C^*$-algebra $A$ so also for the action
($A\ni a \mapsto UaU^*$ and $A\ni a \mapsto U^*aU$) up to the
obtention of a {\em commutative} coefficient algebra $\cal A$ and the
corresponding action (in fact the partial action) which leads to the
corresponding crossed product (in fact the partial crossed product).
The Kwasniewski's crossed product should be naturally considered as
the most general crossed product of the type presented in Example
\ref{1}.

We emphasize that there is {\em no} transfer operator among the
starting objects of this construction and therefore one should
consider the Kwasniewski's construction to be {\em qualitatively}
different from the Exel's crossed product.

\medbreak
Hereafter we describe in brief the Kwasniewski's construction and its
interrelation with the crossed product introduced in the present
article.

Let $A$ be a commutative unital $C^*$-algebra and $\delta$ be an
endomorphism of $A$. As $A$ is commutative we can use the Gelfand
transform in order to identify $A$ with the algebra $C(X)$ of
continuous functions on the maximal ideal space $X$ of $A$. Within
this identification the endomorphism $\delta$ generates (see, for
example, \cite[Proposition 2.1]{Kwa-Leb}) a continuous mapping
$\gamma\!:\Delta \rightarrow X$ where $\Delta\subset X$ is closed and
open (briefly clopen) and the following formula holds
\begin{equation} \label{e3.0}
\delta (a)(x)=
\begin{cases}
a(\gamma(x)),&  x\in \Delta,\\[3pt]
0, &  x\notin \Delta,
\end{cases}
\qquad a\in C(X).
\end{equation}
The mappings $\gamma$ of this sort are called {\em partial mappings}
(of $X$). Thus we have the one-to-one correspondence between the pairs
$(A,\delta)$ and the pairs $(X,\gamma)$, where $X$ is compact and
$\gamma$ is a partial continuous mapping with a clopen domain. In
\cite{kwa} $(X,\gamma)$ is called a \emph{partial dynamical system}.

In \cite{kwa} the author developed the crossed product construction
for an 'almost arbitrary' endomorphism $\delta$, namely the only
presumed constraint was that the {\em image $\gamma (\Delta )$ of the
partial mapping $\gamma$ is open.} As Proposition \ref{prop(O)sition}
and the note preceding it tell this constraint is absolutely minor.

If $A$ is a unital commutative $C^*$-subalgebra of $L(H)$ for some
Hilbert space $H$ and $U$ is a partial isometry such that $U^*U\in A'$
and $UA U^*\subset A$ then clearly $\delta(\,\cdot\,)=U(\,\cdot\,)U^*$
is an endomorphism of $A$. Though in this situation $A$ is not
necessarily a coefficient algebra of $C^* (A,U)$ there is a natural
way to construct one by passage to a bigger $C^*$-algebra ${\cal A}$
generated by $\{A,\,U^*A U,\,U^{2*}A U^2,\,\dotsc\}$. This is stated
in the next

\begin{Pn}\label{prop1} {\upshape \cite[Proposition 4.1]{Leb-Odz}}\ \
If $\delta(\,\cdot\,)=U(\,\cdot\,)U^*$ is an endomorphism of $A$,
$U^*U\in A'$ ($A'$ is the commutant of $A$), and ${\cal
A}=C^*(\bigcup_{n=0}^\infty \delta_*^n(A))$ is the $C^*$-algebra
generated by $\bigcup_{n=0}^\infty U^{*n} A U^n$, then ${\cal A}$ is
commutative and both the mappings \,$\delta\!:{\cal A}\rightarrow
{\cal A}$ and \,$\delta_*\!:{\cal A}\rightarrow {\cal A}$\ \
$(\delta_* (\,\cdot\,) = U^* (\,\cdot\,)U)$ are endomorphisms.
\end{Pn}

\noindent
Thus in the situation described ${\cal A}$ is the
coefficient algebra for  $C^* ({\cal A}, U)= C^*(A,U)$.

It is of primary importance here that ${\cal A}$ is {\em commutative}
and that due to \cite{Kwa-Leb} we can describe its maximal ideal
space, denoted further by $\cal X$, in terms of the maximal ideals in
$A$. Let us recall this description.

To start  with we have to introduce some notation.

Hereafter in this subsection $A$ denotes a commutative unital
$C^*$-algebra, $X$ denotes its maximal ideal space (i.\,e.\ a compact
topological space), $\delta$ is an endomorphism of $A$, while $\gamma$
stands for the continuous partial mapping $\gamma \!:\Delta
\rightarrow X$ where $\Delta\subset X$ is clopen and formula
\eqref{e3.0} holds. When dealing with the partial mappings $\gamma^n$,
\,$n=0,\,1,\,2,\,\dotsc$ we denote, for $n>0$, the domain of
$\gamma^n$ by $\Delta_n=\gamma^{-n}(X)$ and its image by
$\Delta_{-n}=\gamma^{n}(\Delta_n)$; for $n=0$, we set $\gamma ^0 =
\mathrm{Id}$,\ \ $\Delta_0=X$ and thus, for $n,m\in \mathbb N$, we have
\begin{gather}\label{b-2}
\gamma^n\!:\Delta_n\rightarrow  \Delta_{-n},\\[6pt]
\label{b-3}
\gamma^n (\gamma^m(x)) = \gamma^{n+m}(x),\qquad x\in\Delta_{n+m}.
\end{gather}
Note  that in terms of the multiplicative  functionals on $A$,
$\gamma$ is given by
\begin{gather} \label{e1.0}
x\in \Delta_1\ \Longleftrightarrow \ x(\delta(1))=1,\\[6pt]
\label{e2.0}
\gamma(x)=x\circ \delta,\qquad x \in\Delta_1.
\end{gather}

In \cite{Kwa-Leb} the authors calculated the maximal ideal space $\cal
X$ of ${\cal A}$ in terms of $(A,\delta)$ or better to say --- in
terms of the generated partial dynamical system $(X,\gamma)$. This
description is presented in Theorem \ref{ideals2}.

With every $\widetilde{x}\in \cal X$ we associate a sequence of
functionals $\xi^n_{\widetilde{x}}\!:{ A}\rightarrow \mathbb C$,\ \
$n\in\mathbb N$, defined by the condition
\begin{equation} \label{xi2}
\xi^n_{\widetilde{x}}(a)=\delta_*^n(a)(\widetilde{x}),\qquad
a\in A.
\end{equation}
The sequence $\xi^n_{\widetilde{x}}$ determines $\widetilde{x}$
uniquely because ${\cal A}=C^*( \bigcup_{n=0}^\infty\delta_*^n (A))$.
Since $\delta_*$ is an endomorphism of ${\cal A}$ the functionals
$\xi^n_{\widetilde{x}}$ are linear and multiplicative on $A$. So
either $\xi^n_{\widetilde{x}}\in X$ ($X$ is the spectrum of $A$) or
$\xi^n_{\widetilde{x}} \equiv 0$. It follows then that the mapping
\begin{equation}\label{xi000'}
{\cal X}\ni {\widetilde{x}} \to
\bigl(\xi^0_{\widetilde{x}}, \xi^1_{\widetilde{x}},\dotsc\bigr)\in
\prod_{n=0}^{\infty} (X\cup\{0\})
\end{equation} is an injection and  the following statement is true.

\begin{Tm}\label{ideals2}{\upshape
\cite[Theorems 3.1 and  3.3]{Kwa-Leb}}\ \
Let $\delta(\,\cdot\,)=U(\,\cdot\,)U^*$ be an endomorphism of ${ A}$,
$ U^*U\in A$, and $\gamma \!:\Delta_1\rightarrow X$ be the partial
mapping generated by $\delta$. Then the maximal ideal space $\cal X$
of the algebra ${\cal A}=C^*( \bigcup_{n=0}^\infty\delta_*^n ({ A}))$
is given by
\begin{gather*}
{\cal X}=\bigcup_{N=0}^{\infty}{\cal X}_N\cup {\cal X}_\infty,\\[6pt]
{\cal X}_N=\bigl\{\widetilde{x}=(x_0,x_1,\dotsc,x_N,0,\dotsc) :
x_n\in \Delta_n,\ \gamma (x_{n})=x_{n-1},\ 1\leq n\leq N,\
x_N\notin \Delta_{-1}\bigr\},\\[6pt]
{\cal X}_\infty=\bigl\{\widetilde{x}=(x_0,x_1,\dotsc) :
x_n\in \Delta_n,\ \gamma (x_{n})=x_{n-1},\ 1\leq n\bigr\}.
\end{gather*}
The topology on $\bigcup_{N=0}^{\infty}{\cal X}_N\cup {\cal
X}_\infty$ is defined  by a fundamental system of neighborhoods of
points $\widetilde{x}\in {\cal X}_N$
\[
O(a_1,\dotsc,a_k,\varepsilon)=\{\widetilde{y}\in {\cal X}_N:
|a_i(x_N)-a_i(y_N)|<\varepsilon, \ i=1,\dotsc,k\}
\]
and of points $\widetilde{x}\in {\cal X}_\infty$ respectively
\[
O(a_1,\dotsc,a_k,n, \varepsilon)=\biggl\{\widetilde{y}\in
\bigcup_{N=n}^{\infty}{\cal X}_N\cup {\cal X}_\infty\,:\,
|a_i(x_n)-a_i(y_n)|<\varepsilon, \ i=1,\dotsc,k\biggr\}
\]
where $\varepsilon > 0$, $a_i\in A$ and $k,n \in \mathbb N$.
\end{Tm}

\begin{Rk}\label{remark2.3}
The topology on $X$ is $^*$-weak. One immediately sees then (see
\eqref{xi000'}), that the topology on $\bigcup_{N\in \mathbb N}{\cal
X}_N\cup {\cal X}_\infty$ is in fact the product topology inherited
from $\prod_{n=0}^{\infty} (X\cup\{0\})$ where $\{0\}$ is clopen.
\end{Rk}

This theorem motivates us to take a closer look at the condition
$U^*U\in A$.

Observe \cite[Proposition 3.5]{Kwa-Leb} that if $U^*U\in A'$ then
$\delta$ is an endomorphism of the $C^*$-algebra ${ A}_1=C^*({ A},
U^*U)$ and we also have
\[
{\cal A}=C^*\left(\bigcup_{n=0}^\infty \delta_*^n(A)\right) =
C^*\left(\bigcup_{n=0}^\infty \delta_*^n(A_1)\right).
\]
Thus the mentioned condition simply means that when calculating the
maximal ideal space $\cal X$ of $\cal A$ one should start from the
$C^*$-algebra $A_1$ rather than from $A$.

Furthermore, the condition $U^*U\in A$ is closely related to the
openness of $\Delta_{-1}$ (as $\Delta_1$ is compact and $\gamma$ is
continuous $\Delta_{-1}$ is always closed).

\begin{Pn}\label{prop(O)sition}
{\upshape \cite[Proposition 2.4]{kwa}}\ \
Let $P_{\Delta_{-1}}\in A$ be the projection corresponding to
characteristic function $\chi_{\Delta_{-1}}\in C(X)$. If  $U^*U\in
A$ then  $\Delta_{-1}$ is open and  $U^*U = P_{\Delta_{-1}}$. If
$U^*U\in A'$,  $\Delta_{-1}$ is open and $A$ acts
non-degenerately on $H$, then $ U^*U \leqslant P_{\Delta_{-1}}.$
\end{Pn}

\noindent
\textbf{Note.} As it is shown in \cite{kwa} the inequality in the
second part of the preceding proposition can not be replaced by
the equality.

By virtue of Proposition \ref{prop1} the mappings $\delta$ and
$\delta_*$ are endomorphisms of the $C^*$-algebra ${\cal A}$. Making a
start from Theorem \ref{ideals2} we can now find the form of the
partial mappings they generate. The initial hint here is the following

\begin{Pn}
\label{tu1} {\upshape \cite[Proposition 2.5]{Kwa-Leb}}\ \
Let $\delta(\,\cdot\,)=U(\,\cdot\,) U^*$, and
$\delta_*(\,\cdot\,)=U^*(\,\cdot\,)U$ be endomorphisms of $A$ and let
$\gamma $ be the partial mapping of $X$ generated by $\delta$. Then
$\Delta_1$ and $\Delta_{-1}$ are clopen and $\gamma\!:\Delta_1
\rightarrow \Delta_{-1}$ is a homeomorphism. Moreover the endomorphism
$\delta_*$ is given on $C(X)$ by the formula
\begin{equation}
(\delta_* f)(x)=
\begin{cases}
f(\gamma^{-1}(x)),&  x\in \Delta_{-1},\\[3pt]
0, &  x\notin \Delta_{-1}.
\end{cases}
\end{equation}
\end{Pn}

This proposition along with Theorem \ref{ideals2} leads to

\begin{Tm}\label{tu3}
{\upshape \cite[Theorem 2.8]{kwa}}\ \
Let the hypotheses of Theorem \ref{ideals2} hold. Then
\begin{itemize}
\item[$i)$] the sets
\begin{gather*}
\widetilde{\Delta}_1=\{\widetilde{x}=(x_0,x_1,\dotsc)\in {\cal X}: x_0
\in \Delta_1\},\\[3pt]
\widetilde{\Delta}_{-1}=\{\widetilde{x}=(x_0,x_1,\dotsc)\in {\cal X}:
x_1\neq 0 \}
\end{gather*}
are clopen subsets of ${\cal X}$,
\item[$ii)$] the endomorphism $\delta$ generates on
${\cal X}$ the partial homeomorphism
$\widetilde{\gamma}\!:\widetilde{\Delta}_1 \rightarrow
\widetilde{\Delta}_{-1}$  given by the formula
\begin{equation}\label{alfazfalkom}
\widetilde{\gamma}(\widetilde{x})=\widetilde{\gamma}(x_0,x_1,\dotsc)
=(\gamma(x_0),x_0,x_1,\dotsc),\qquad
\widetilde{x}\in \widetilde{\Delta}_1,
\end{equation}
\item[$iii)$]
the partial mapping generated by  $\delta_*$ is the inverse of
$\widetilde{\gamma}$, that is
$\widetilde{\gamma}^{-1}\!:\widetilde{\Delta}_{-1}\rightarrow
\widetilde{\Delta}_{1}$ where
\begin{equation}\label{alfa-1zfal}
\widetilde{\gamma}^{-1}(\widetilde{x})=\widetilde{\gamma}^{-1}
(x_0,x_1,\dotsc)=(x_1,x_2,\dotsc),\qquad
\widetilde{x}\in \widetilde{\Delta}_{-1}. \end{equation}
\end{itemize}
\end{Tm}

The pair $({\cal X}, \widetilde{\gamma})$ is called in \cite{kwa} the
{\em reversible extension of the partial dynamical system}
$(X,\gamma)$.

Along with the description of the maximal ideal space of $\cal A
=C({\cal X})$ and the action on it the author presented in \cite{kwa}
an explicit algebraic construction of $\cal A$ in terms of
$(A,\delta)$. Here is this construction.

Observe first that the family $\{A_n\}_{n\in\mathbb N}$ where
$A_n:=\delta^n(1)A$, $n\in \mathbb N$, is a decreasing family, of
closed two-sided ideals. Since the operator $\delta^n(1)$ corresponds
to $\chi_{\Delta_n}\in C(X)$, one can consider $A_n$ as
$C_{\Delta_n}(X)$ (we denote by $C_{K}(X)$ the algebra of continuous
functions on $X$ vanishing outside the set $K\subset X$). Let $
\mathcal{E}_*(A)$ be the set consisting of the sequences
$a=\{a_n\}_{n\in \mathbb N}$ where $a_n\in A_n$, $n=0,\,1,\,\dotsc$,
and only a finite number of functions $a_n$ are non zero. Namely
\[
\mathcal{E}_*(A)=\biggl\{a\in \prod_{n=0}^\infty A_n: \exists\,{N>0}\
\forall\,{n>N}\quad a_n\equiv 0\biggr\}.
\]
Let $a=\{a_n\}_{n\geqslant 0}$, $b=\{b_n\}_{n\geqslant 0}\in
\mathcal{E}_*(A)$ and $\lambda\in\mathbb{C}$. We define the
addition, multiplication by scalar, convolution multiplication and
involution on $ \mathcal{E}_*(A)$  as follows
\begin{gather}\label{add}
 (a+b)_n=a_n+b_n,\\[6pt]
\label{mulscal}
 (\lambda a)_n=\lambda a_n,\\[6pt]
\label{mul}
 (a\,\cdot\, b)_n=a_n\sum_{j=0}^n \delta^j(b_{n-j})+b_n \sum_{j=1}^n
\delta^j(a_{n-j}),\\[6pt]
\label{invol}
 (a^*)_n=\overline{a}_n.
\end{gather}
These operations are well defined and very natural, except maybe the
multiplication of two elements from $\mathcal{E}_*(A)$. We point out
here that the index in one of the sums of \eqref{mul} starts running
from $0$.

\begin{Pn}\label{algebrawithinvolution}
{\upshape \cite[Proposition 4.5]{kwa}}\ \
The set $\mathcal{E}_*(A)$ with operations \eqref{add},
\eqref{mulscal}, \eqref{mul}, \eqref{invol} becomes a commutative
algebra with involution.
\end{Pn}

Now, let us define a morphism $\varphi\!:\mathcal{E}_*(A)\rightarrow
C({\cal X})$. To this end, let $a=\{a_n\}_{n\in \mathbb N}\in
\mathcal{E}_*(A)$ and $\tilde{x}=(x_0,x_1,\dotsc)\in \cal X$. We set
\begin{equation}\label{phimap}
\varphi(a)(\tilde{x})=\sum_{n=0}^\infty a_n(x_n),
\end{equation}
where $a_n(x_n)=0$ whenever $x_n=0$. The mapping $\varphi$ is well
defined as only a finite number of functions $a_n$, $n\in \mathbb N$,
are non zero.

\begin{Tm}\label{E*} {\upshape \cite[Theorem  4.6]{kwa}}\ \
The mapping $\varphi\!:\mathcal{E}_*(A)\rightarrow C({\cal X})$ given
by \eqref{phimap} is a morphism of the algebras with involution. The
image of $\varphi$ is dense in $C({\cal X})$, that is
\[
\overline{\varphi(\mathcal{E}_*(A))}=C({\cal X}).
\]
\end{Tm}

Let us consider the quotient space $\mathcal{E}_*(A)/
\mathop{\mathrm{Ker}}\varphi$ and the corresponding quotient mapping
$\phi\!:\mathcal{E}_*(A)/\mathop{\mathrm{Ker}}\varphi\rightarrow
C({\cal X})$, that is $\phi(a+\mathop{\mathrm{Ker}} \varphi)
=\varphi(a)$. Clearly $\phi$ is an injective mapping onto a dense
$^*$-subalgebra of $C({\cal X})$. Let us set
\begin{gather}\label{e(c)}
E_*(A):=\phi(\mathcal{E}_*(A)/\mathop{\mathrm{Ker}}\varphi),\\[6pt]
\label{phi}
[a]:=\phi(a+\mathrm{Ker}\,\varphi), \qquad a\in \mathcal{E}_*(A)
\end{gather}

\begin{Dn}\label{idenitify-defn}
{\upshape \cite[Definition 4.7]{kwa}}\ \
$E_*(A)$ is called a \emph{coefficient\/ *-algebra} of the
pair $(A,\delta)$. We shall write $[a]=[a_0,a_1,\dotsc]\in
C({\cal X})$ for $a=(a_0,a_1,a_2,\dotsc)\in \mathcal{E}_*(A)$.
\end{Dn}

The natural injection $ A\ni a_0\longrightarrow [a_0,0,0,\dotsc]\in
E_*(A) $ enables us to treat $A$ as a $C^*$-subalgebra of $E_*(A)$ and
hence also of ${\cal A} = C({\cal X})$:
\[
A\subset E_*(A)\subset {\cal A},\qquad \overline{E_*(A)}={\cal A}.
\]

Once the extension of $A$ up to the coefficient algebra $\cal A$
is implemented  the further construction of the crossed product in
\cite[Section 5]{kwa} goes smoothly. Since here $\delta$ and
$\delta_*$ are the endomorphisms of $\cal A$ they are partial
automorphisms (Proposition \ref{tu1}) and their actions are
described in Theorem \ref{tu3}. Therefore the crossed product in
\cite{kwa} is naturally defined as the partial crossed product
developed in \cite{Exel1}. Clearly in the situation considered it
coincides with ${\cal A}\times_{\delta} {\mathbb Z}$.

\medbreak
To finish this subsection we would like to mention a number of
interesting observations made in \cite{kwa} about the dependence
of $\cal X$  on $\gamma$.

If $\gamma$ is surjective then ${\cal X}_N$, $n\in \mathbb N$, are
empty and ${\cal X}={\cal X}_\infty$, in this case $\cal X$ can be
defined as a projective limit \cite[Proposition 3.10]{kwa}.

If $\gamma$ is injective then a natural continuous projection $\Phi$
of $\cal X$ onto $X$ given by the formula
\begin{equation} \label{Phimap} \Phi(x_0,x_1,\dotsc) =x_0
\end{equation}
is a homeomorphism \cite[Proposition 2.3]{kwa}.

If $(X, \gamma)$ is the {\em one-sided} topological Markov chain then
its reversible extension $({\cal X}, \tilde{\gamma})$ is the {\em
two-sided} topological Markov chain \cite[Example 2.8]{kwa}.

The latter example shows in particular that the Kwasniewski's crossed
product is qualitatively different from the Cuntz-Krieger algebra.
Though both these crossed products start from the one-sided
topological Markov chain in the Cuntz-Krieger construction we arrive
at the coefficient algebra ${\cal F}_A$ which is not commutative (see
Subsection \ref{C-K}) while the Kwasniewski's construction leads to
the commutative coefficient algebra $C({\cal X})$ where $({\cal X},
\tilde{\gamma})$ is the two-sided topological Markov chain.

\section{Crossed product, coefficient algebras,\\
transfer operators, etc (interrelations) }
\label{crossed product etc}

The constructions presented in the previous section being different
(and sometimes even 'qualitatively' different, here we would like to
'counterpose' Exel's and Kwasniewski's crossed products) all turned
out to be of the type of the crossed product introduced in this
article: in all of them there arise a finely representable pair $(\cal
A,\delta)$ and the corresponding copmlete transfer operator~$\delta_*$.

But since nevertheless these constructions are {\em different\/} one
naturally arrives at the question: what object should be named the
'crossed product'? Or may be 'more exact' question: 'what' are we
crossing with 'what'?

The objects of the previous section along with the results of Sections
\ref{crossed product} and \ref{faithful} lead to the 'natural answer'
which we are giving hereafter.

\begin{I}\label{cr-pr-scheme}
{\bfseries Crossed product construction.}\ \
The crossed product construction consists of {\em two\/} steps.

\smallbreak
{\bfseries Step 1 (Initial object and extension procedure).}\ \
There should be given a certain $^*$-algebra $A$ (or even only a
certain set of elements of $A$) and an {\em extension procedure} by
means of which one can extend the algebra $A$ up to a {\em
coefficient\/} $C^*$-algebra $\cal A$ and define an endomorphism
$\delta \!: {\cal A} \to {\cal A}$ in such a way that the pair
$(\cal A,\delta)$ is finely representable (this is equivalent to
the existence of a complete transfer operator $\delta_*$ for
$({\cal A},\delta)$ which is unique by Theorem~\ref{complete}).

\smallbreak
{\bfseries Step 2 (Crossing the coefficient algebra with the
endomorphism $\delta$).}\ \
Once ${\cal A}$, $\delta$ and $\delta_*$ are given the crossed product
(of $\cal A$ and $\delta$) is defined according to
Definition~\ref{crossed}.
\end{I}

\medbreak
So 'in essence' we are crossing the {\em coefficient algebra $\cal A$}
with the {\em endomorphism} $\delta$, while the initial algebra $A$
and the extension procedure serve as an instrument to construct $\cal
A$ and $\delta $.

\begin{I}\label{discus}{\bfseries Discussion.}

\medbreak
{\bfseries I\ \,Extension.}\ \
In step 1 of the crossed product construction the extension procedure
was mentioned. What is it? Is there any possibility to describe it
explicitly? In the general setting (unfortunately) we do not see the
way to give the complete description. Looking through the crossed
product structures considered in Section \ref{v-crossed product} one
can notice that neither a transfer operator (see, for example, the
Kwasniewski's crossed product) nor even an endomorphism of $A$ (see
Example \ref{22}) are necessary for this procedure. Moreover if we
look at the Cuntz-Krieger algebra then among the starting objects we
find the projections $Q_i =S^*_iS_i$ and $P_i = S_iS^*_i$ the set of
which {\em does not} even form an {\em algebra} and there is {\em no
transfer operator} among the starting objects as well (the necessary
algebra and transfer operator arise here as the outcome of the
extension procedure). We would like also to mention in this connection
the paper \cite{Lin-Rae} by J. Lindiarni and I. Raeburn where the
starting objects are a $C^*$-algebra $A$ and a positive cone
$\Gamma^+$ of a totally ordered abelian group acting on $A$ by
endomorphisms (that are extendible onto the multiplier algebra of
$A$); {\em no transfer operator(s)} is(are) presumed. Assuming that
these endomorphisms are generated by partial isometries $V_s$,
\,$s\in \Gamma^+$ in such a way that the endomorphisms are given by
$A\ni a \mapsto V_saV^*_s$ and $V^*_sV_s\in A^\prime$ the authors
define the crossed product as the universal enveloping $C^*$-algebra
generated by $A$ and $V_s$. Clearly here the corresponding coefficient
algebra $\cal A$ (the extension of $A$) is the universal $C^*$-algebra
generated by $A$, $V^*_sAV_s$, \ $s\in \Gamma^+$ and $\delta_s
(\,\cdot\,) = V_s (\,\cdot\,)V^*_s$ and $\delta_{*s} (\,\cdot\,) =
V^*_s (\,\cdot\,)V_s$.

In fact in the general situation the extension procedure is given by
the mappings $A\ni a \mapsto UaU^*$ and $A \ni a \mapsto U^*aU$ that
{\em depend} on the {\em origin} of the algebra $A$ and the {\em
origin} of the operator $U$ (in general neither the mapping $U
(\,\cdot\,)U^*$ nor the mapping $U^*(\,\cdot\,)U$ preserve $A$ so they
are neither endomorphisms nor transfer operators). In principle the
{\em extension procedure} reduces to a certain {\em axiomatic
description of these mappings}. The main structures and objects that
will appear on this way are described in \cite[Section 3]{Leb-Odz}.

Here we would like also to remind the already mentioned 'contrast'
between the Cuntz-Krieger algebra and the Kwasniewski's crossed
product for the one-sided topological Markov chain. Both these crossed
products start from the one-sided topological Markov chain. Then the
operators $S_i$ in the Cuntz-Krieger algebra generate the isometry $S$
\eqref{S} that leads to the extension procedure given by
$S(\,\cdot\,)S^*$ and $S^* (\,\cdot\,)S$ and ending with the
coefficient algebra ${\cal F}_A$; while the extension procedure for
the Kwasniewski's crossed product leads to the commutative coefficient
algebra $C({\cal X})$ where $({\cal X}, \tilde{\gamma})$ is the
two-sided topological Markov chain.

In the recent paper \cite{exel-inter} by R.\ Exel certain extension
procedures are discussed that have the form
\[
U^* a U = {\cal H} (a)U^*U \quad \text{and}\quad U a U^* = {\cal V}
(a)UU^* , \qquad a\in A
\]
where $A$ is a $C^*$-algebra, $U$ is a partial isometry and ${\cal V,
H}\!: A \to A$ are some positive linear maps. After finishing a series
of fascinating calculation at the end of \cite[Section 7]{exel-inter}
R.\ Exel exclaims: 'The reader may be struck with the impression that
the wild juggling of covariant representations ... is a bit
exaggerated and that something must be done to stop it. I agree. ...'
We agree as well. But all the foregoing reasoning convinces us that
there is {\em no} universal way to describe an {\em arbitrary}
extension procedure. We repeat once more that there are as many
extension procedures as types of $C^*$-algebras $A$ and mappings $A\ni
a \mapsto UaU^*$ and $A \ni a \mapsto U^*aU$.

\medbreak
{\bfseries II\ \,Coefficient algebras and transfer operators.}\ \
In this article we have {\em started\/} with the coefficient algebras
(Section \ref{crossed product}). But in fact as it is clear from the
above discussion the coefficient algebra is not a starting object but
rather the {\em ending} (intermediate) one --- it is the {\em
target\/} of the extension procedure.

Concerning the interrelations between the coefficient algebras and
transfer operators we have to emphasize that {\em at the beginning\/}
(when we are starting to construct a coefficient algebra) a transfer
operator {\em may or may not\/} arise (it depends on the extension
procedure we are confronted with) while {\em at the end\/} (when we
have already constructed the coefficient algebra) this operator
appears {\em necessarily\/} and is {\em unique}.

It is worth mentioning that in the Exel's construction of the crossed
product (see \ref{Ex}) the transfer operator $\cal L$ is {\em not\/}
defined in a {\em unique\/} way (it is not totally defined by
$\alpha$, see, for example, \cite[Remark 2.4]{Bakht-Leb}). On the
other hand when by means of the Exel's construction one obtains a
coefficient algebra $\cal A$ (like in Theorems \ref{Tau-cr} and
\ref{Ex!}) then the arising transfer operator $\delta_*$ is {\em
unique\/} (it is totally defined by $\delta$ (Theorem
\ref{complete})). Therefore the difference between various Exel's
transfer operators $\cal L$ indicates the {\em difference\/} between
the {\em starting\/} objects for the corresponding various {\em
extension\/} procedures (we repeat that {\em not all\/} of the
extension procedures can be defined by means of a certain transfer
operator).

\medbreak
{\bfseries III\ \,Crossed product.}\ \
The construction given in Section \ref{crossed product} {\em does
not\/} cover all the already existing crossed product type
constructions (in particular the Exel's crossed product associated
with the transfer operator $\cal L$ that is not related to a partial
isometry is not embedded into it). On the other hand we believe (and
the material of the article convinces us in this) that once one is
confronted with the crossed product structure related to {\em partial
isometries\/} he will necessarily come to the crossed product
construction described in \ref{cr-pr-scheme} and therefore at the end
to the crossed product given in Definition \ref{crossed}.

This crossed product is quite satisfactory in a number of ways: it
covers all the most successful crossed product structures
developed earlier (see Section \ref{v-crossed product}), and it
possesses good internal structural properties (described in
Section \ref{faithful}).

Therefore recalling Exel's exclamation: '...something must be
done to stop it...' we beleive that this crossed product may be
considered as a reasonable (at least intermediate) stop.
\end{I}

\end{document}